\documentclass[12pt,reqno,letter]{amsart}
\newif\iflong
\longtrue
\usepackage{amsmath,amsfonts,amssymb,amscd,amsthm,amsbsy}
\usepackage[inline]{enumitem}

\usepackage[utf8]{inputenc}
\usepackage{caption}
\usepackage{graphicx}
\graphicspath{{images/}}

\usepackage{hyperref}
\usepackage{multirow}
\usepackage[dvips]{xcolor}

\usepackage{longtable}
\allowdisplaybreaks

\title[Break-down of quasi-periodic attractors] {Accurate computations
up to break-down of quasi-periodic attractors in the dissipative
spin-orbit problem}

\author[R. Calleja]{Renato Calleja}

\address{ Department of Mathematics and Mechanics, IIMAS, National
  Autonomous University of Mexico (UNAM), Apdo. Postal 20-126, Mexico
  D.F., Mexico \& Departament de Matem\`atiques i Inform\`atica,
  Universitat de Barcelona, Gran Via de les Corts Catalanes, 585,
  08007 Barcelona, Spain}

\email{celleja@mym.iimas.unam.mx}

\author[A. Celletti]{Alessandra Celletti}

\address{ Department of Mathematics, University of Rome Tor Vergata,
  Via della Ricerca Scientifica 1, 00133 Rome, Italy}

\email{celletti@mat.uniroma2.it}

\author[J. Gimeno]{Joan Gimeno}

\address{ Departament de Matem\`atiques i Inform\`atica, Universitat
  de Barcelona, Gran Via de les Corts Catalanes, 585, 08007 Barcelona,
  Spain}

\email{jgimeno@ub.edu}

\author[R. de la Llave]{Rafael de la Llave}

\address{ School of Mathematics, Georgia Institute of Technology, 686
  Cherry St., Atlanta, GA 30332-0160, USA }

\email{rafael.delallave@math.gatech.edu}

\thanks{R.C. acknowledges the hospitality of the Department of
  Mathematics at the University of Barcelona, where he had a Mar\'ia
  Zambrano Fellowship. R.C. was also partially supported by DGAPA-UNAM project PAPIIT IN103423 and PASPA. A.C. was partially supported by EU H2020 MSCA
  ETN Stardust-Reloaded Grant Agreement 813644. J.G. has been
  supported by the Italian grant MIUR-PRIN 20178CJA2B ``New Frontiers
  of Celestial Mechanics: theory and Applications'', with funds
  from NextGenerationEU within the Spanish national Recovery,
  Transformation and Resilience plan, with the Spanish
grant PID2021-125535NB-I00 (MICINN/AEI/FEDER, UE), and the project that gave rise to these results received the support of a fellowship from ``la
Caixa'' Foundation (ID 100010434). The fellowship code is 129200.  J.G. also thanks the School of
  Mathematics of GT for its hospitality in Fall 2019 and year
  2022. R.L. has been partially supported by NSF grant DMS 1800241.}

\date{\today}

\newcommand{\I}{\mathtt{i}}

\newcommand{\eps}{\varepsilon}
\newcommand{\dis}{\mu}
\newcommand{\ecc}{e}
\newcommand{\omg}{\omega}

\newtheorem{meta-thm}[thm]{Meta-Theorem}

\newtheorem{defn}{Definition}

\newtheorem{rem}{Remark}

\newtheorem{alg}{Algorithm}

\newcommand{\bydef}{\,\stackrel{\mbox{\tiny\textnormal{\raisebox{0ex}[0ex][0ex]{def}}}}{=}\,}

\newcommand\beq[1]{ \begin{equation}\label{#1} }
\newcommand{\eeq}{ \end{equation} }

\newcommand\beqa[1]{ \begin{eqnarray} \label{#1}}
\newcommand{\eeqa}{ \end{eqnarray} }
\newcommand{\beqano}{ \begin{eqnarray*} }
\newcommand{\eeqano}{ \end{eqnarray*} }
\newcommand\equ[1]{{\rm (\ref{#1})}}
\newcommand{\R}{\mathbb{R}}
\newcommand{\Z}{\mathbb{Z}}
\newcommand{\N}{\mathbb{N}}
\newcommand{\T}{\mathbb{T}}

\DeclareMathOperator{\Tr}{Tr}

\begin{document}


\maketitle \vspace{-1em}
\tableofcontents

\begin{abstract}
We consider a Celestial Mechanics model: the spin-orbit problem
with a dissipative tidal torque, which is a singular perturbation
of a conservative system.  The goal of this paper is to show that
it is possible to compute quasi-periodic attractors accurately and
reliably  for parameter values extremely close to the breakdown.
Therefore, it is possible to obtain information on mathematical
phenomena at breakdown.

The method we use incorporates  the same time numerical and rigorous improvements.
Among them
\begin{enumerate*}[label={\sl (\roman*)}]
 \item the formalism is  based on studying the
  time-one map of the spin-orbit problem (which reduces the dimensionality of the problem)
and has mathematical advantages;
\item very accurate integration of the ODE (high order Taylor methods implemented with extended precision)
for the map at its jets;
 \item a very efficient KAM method for maps which computes the attractor and its tangent spaces ( quadratically
   convergent step with  low storage requirements, and low operation
   count);
 \item the algorithms are backed by a  rigorous a-posteriori KAM
   Theorem, which establishes that if the algorithm,
   produces a very approximate solution of functional equation with reasonable
condition numbers.   then there is a true solution nearby;
   and
 \item   the continuation algorithm is  guaranteed to reach arbitrarily close to the
   border of existence if it is given enough computer resources.
\end{enumerate*}

As a byproduct of the accuracy that we maintain till breakdown, we study
several scale invariant observables of the tori used in the renormalization group of infinite dimensional spaces.  In contrast with
previously studied simple models, the behavior at breakdown of the
spin-orbit problem does not satisfy standard scaling relations which implies that the spin-orbit problem is not
described by a hyperbolic fixed point of a  renormalization
operator.
\end{abstract}

\section{Introduction}
The existence of invariant tori of maximal dimension is important in
the study of the stability of dynamical systems; when the number of
degrees of freedom is low, invariant tori act as barriers and confine
the motion in regions of phase space.

In recent papers (\cite{CCGL20a,CCGL20b}), we have developed extremely
effective methods to compute very accurately quasi-periodic solutions
(KAM tori) of several models of the dissipative spin-orbit problem in
celestial mechanics.

The spin-orbit problem describes the motion of an oblate satellite
whose center of mass moves on a Keplerian orbit around a central
planet. The satellite is not rigid and its motion generates tidal
friction. We consider both a model in which the friction is time
dependent along the position in the orbit (the tides are different at
different places in the orbit) as well as another one whose friction
is constant (an average of the time dependent friction).

Both models depend on three parameters: the oblateness of the
satellite, the orbital eccentricity, and the dissipative factor.  The
limit of zero dissipation is a singular limit, since the long term
behaviors of the model change drastically from a conservative system
to a dissipative one even if the dissipation is small.

For dissipative systems, to obtain a quasi-periodic orbit of a certain
frequency, one has to adjust a parameter called drift. This is a big contrast with the conservative case which
admits many
quasi-periodic orbits. To obtain orbits with a given frequency in a conservative system, one
just needs to choose the initial conditions and one does not need to adjust parameters.

The astronomically relevant values of the friction are typically
small, so that the systems of astronomical interest are in the regime of
singular perturbations (frictions have crucial effects in the
formation of planets over millions of years, but are very small in
calculations over a few thousand years).
Computing by direct iterations  requires long transients.

Note that numerical methods based on contractions (e.g. the numerical
methods suggested in \cite{CCFL14}) slow down in the limit of
small dissipation, whereas the Newton methods we use maintain the speed in
the limit of small dissipation.

Working close to a singular limit, our numerical explorations,
especially when we want to study the objects near the breakdown,
require several elaborated numerical methods (see
Sections~\ref{sec:rotation} and \ref{sec:breakdown}) as well as some
theoretical underpinnings, e.g. condition numbers and theoretical
arguments that indicate that the numerical calculations are
credible. Both the efficient numerical algorithms and the a-posteriori
results take advantage of several identities that are obtained using
the geometric properties of a (conformally symplectic) map.

The goal of this paper is to show empirically that our methods are
able to continue very effectively the quasi-periodic solutions
w.r.t. the parameters; at the same time, our methods allow to maintain
a very high accuracy, even extremely close to parameter values where
the solutions cease to exist. Our work gives results in different
directions:
\begin{enumerate*}[label={\sl (\roman*)}]
 \item It provides a very good test of the methodology to maintain
   accuracy and reliability;
 \item It shows that the \emph{a-posteriori} theorems developed in
   \cite{CallejaCdlL2013} can be used in very practical situations and
   effectively cooperate with elaborated numerical methods;
 \item It yields information of mathematical phenomena happening at
   breakdown of the invariant tori.
\end{enumerate*}

One of our main tools to ascertain the breakdown is that the Sobolev
seminorms of the embedding functions describing the tori blow up when
the parameters approach the breakdown.  We also study the behavior of
the stable and tangent bundles; in \cite{CallejaF2012}, it was found
numerically that at the breakdown of tori in the dissipative standard
map, one has simultaneously a loss of regularity of the torus and a
bundle collapse between the stable and tangent bundles (the center and the stable bundles get close in a dense set, even if they have to remain
apart in a set of full measure, hence they have to become very
oscillatory).  We note that there are rigorous results
\cite{CallejaL10, CallejaCdlL2013} which show that the algorithms can
compute arbitrarily close to breakdown and that the breakdown is
signaled either by a blow up of Sobolev seminorms of the conjugacy or
a loss of hyperbolicity (in conformally symplectic systems, the stable
Lyapunov exponent is bounded away from zero, so that the loss of
hyperbolicity can only happen because the minimum of the angles between the stable
and tangent bundles are not bounded away from zero).

In contrast to previous observations in other models (mainly maps with
very few harmonics), we have found that the Sobolev seminorms do not
have an asymptotic power law (with log-periodic corrections)  and that the bundles do not oscillate
wildly (they remain separated even if they lose regularity). There do
not seem to occur standard scaling relations on the properties of the
bundles.

The failure of the scaling behaviour found in other
models with few harmonics
can  be described in the  renormalization group  language as saying
that the spin-orbit problem is so far away from previously studied models
that it is in
a different universality class.   Indeed, \cite{MacKay82}
asked the question of clarifying the limits of the universality class.

It seems quite likely that some of the phenomena discovered here for
the spin-orbit problem should be applicable to more models (they are
explained by certain properties of the dynamics of the renormalization
group). Since the goal of this paper is to show how one can maintain
accuracy and reliability up to breakdown in a concrete model of
astronomical interest, we have not done a systematic exploration of
similar models and we postpone the worthwhile goal of understanding dynamics of the renormalization group.

\subsection{Some comments on the numerical tools}

To reach our results at the limit of breakdown, we have deployed the
following numerical tools: a reduction of the spin-orbit problem to a
time 1-map, a very high order numerical integrator, jet transport to
compute the derivatives, and extended precision arithmetic.  Let us
briefly comment on these tools.

\subsubsection{Reduction to time-1 map}
The reduction to a time-1 map encodes the KAM tori by a lower
dimensional object. The computation of the time-1 map decreases the
memory usage but, in our case, it requires several integrations of an
ODE, which is computationally expensive. However, the resulting code
is very parallelizable since the integration of each orbit can be
assigned to a thread (we have used this systematically).

Note that the tori --specially near breakdown-- become rough along the
transversal direction while they remain smooth along the flow
direction.  Hence, it is natural to deal with these directions
differently and to treat the integration along the flow (which remains
always very smooth) differently than the KAM computation.

Another advantage of basing the results on the time-1 map is that it
permits a rigorous justification of the comparison between the
spin-orbit model \eqref{eq.spin-diss} and the averaged version
\eqref{eq.avg-spinxy} (see Remark~\ref{averages}).

\subsubsection{Using high order (Taylor) methods, jet transport and extended precision; the {\tt taylor} package}
As it turns out, the KAM method requires very high precision in the
torus which can be achieved only  by using extended precision in the
arithmetic and in the integration step. The only way to achieve
sufficient precision in the integration using reasonable size steps
is to use a very high order method (e.g. order 30).

Moreover, the KAM algorithm uses the  derivatives with respect to initial conditions
and parameters of the time one map,
which requires integrating the variational
equations.

All these requirements (very high precision calculations for the map and
the variational equations)
can be obtained simultaneously with a
reasonable amount of programming time using the public
domain {\tt taylor} package \cite{JorbaZ2005,GimenoJZ22}.

The {\tt taylor} program takes as input the equations written in a
simple symbolic form and it outputs a C/C++ code that implements a
Taylor stepper integrator.  One of the options of {\tt taylor} is to
implement the arithmetic in an extended precision package by using,
for instance, the {\tt MPFR} package. If the user is not satisfied
with the constant step, {\tt taylor} allows to select and fine tuning
several adaptive steps or even implement his/her own time stepping
strategies.

We recall that the Taylor integrator (already used by I. Newton)
consists in finding recursively the Taylor coefficients of the
solution.  The differential equation gives a recursion that prescribes
the Taylor coefficients of the solution of order $n$ as a function of
the coefficients of lower orders.  The initial conditions are the
starting point of these recursions, which are easy to formulate when the
differential equations are formed from elemental functions, such as;
$\sin, \cos, \exp,\log, \text{pow}$, etc.  and using composition and
arithmetic operations.

An important  recently added   feature of {\tt taylor} described in
\cite{GimenoJZ22} is that it implements what is called {\em jet
  transport} and it provides high-order derivatives of the time-$t$
map. Jet transport is a technique that considers the phase space
variables as  multivariate polynomials instead of real
numbers.   As it is proved in \cite{GimenoJJMZ22} the numerical solution
of the  ODE with polynomial variables is equivalent to solving the
original equation and its variational equations. The Taylor method
for polynomials can be obtained by overloading the polynomial
types in the arithmetic operations and one can also adapt
step control strategies.

The upshot is that, using {\tt taylor}, with a programming effort
comparable to typing the equations in \TeX, one obtains an arbitrary
order solver for the equations and the variational equations,
implemented in extended precision. As we will see, the running times
are very competitive.

In our case, we computed the one-time map and its derivatives with a
precision which is (conservatively) smaller than $10^{-35}$ in a
reasonable computer time. Higher precisions could be
accomplished by changing a few parameters in the run at the price of
more computer resources.

\subsubsection{Computation of KAM tori given the map; automatic reducibility}
Once we have the time-1 map very carefully computed, we turn on the
problem of computing the KAM tori by a continuation scheme.

The method is based on formulating a functional equation for
the parameterization of a torus that expresses invariance
and the motion is a  rotation.  The unknowns to be solved for
are the parameterization and the parameters of the model
that need to be adjusted.

We use a Newton-like method based on automatic reducibility developed
for symplectic mappings in \cite{Llave01c, LlaveGJV05} and for
conformally symplectic systems in \cite{CallejaCdlL2013}.

Note that the algorithm only needs to manipulate functions of
as many variables as the dimension of the torus in the
return map.
This number of variables is smaller than the
dimension of the phase space (other methods based on
normal forms or expansions use functions in phase space
or even families of such functions).  Reducing the number of variables of
the functions considered is important because the computational  cost of handling
a function grows extremely rapidly with the number of variables
(the \emph{curse of dimensionality}).

The algorithm we use  takes advantage of several geometric identities to
obtain a quadratically convergent method. The method has very low
storage requirements and operation count (if the torus is discretized
with $N$ points, a step requires $O(N)$ storage and $O(N\log(N))$
operations).  The algorithm (listed in \cite{CallejaCdlL2013})
consists in about a dozen steps on objects. Using modern computer
languages, all of those steps can be implemented in one line in a way
that is independent of the level of truncation or accommodating
extended precision arithmetic. Using more traditional languages,
one can implement a layer of operations among objects
and then a sort main program.  Hence, significant parts of the problem
can be reused.

\subsubsection{A posteriori results}

Since we are going to compute very close to breakdown, reliability is
not obvious.

An important consideration is the rigorous result of \cite{CallejaCdlL2013} stated in an a-posteriori format, that is: given an
approximate solution, if one computes some (explicit) condition numbers,
then one can ensure that there is a true solution nearby. Hence,
provided that we can compute solutions with good accuracy and keeping
track of the condition numbers, we are confident of our calculations
to breakdown.  Explicit (but not completely rigorous since they
ignored round off) calculations of the condition numbers which give
existence of the tori extremely close to the breakdown appear in
\cite{CCGL20b}. This has been reconsiderd in \cite{Linroth23}.

We also note that the fact that we have a-posteriori results of the
right form tells us that a continuation method will progress unless
the torus ceases to satisfy a few conditions. Hence, the method
guarantees that our continuation methods will reach arbitrarily close
to breakdown (if we could use enough discretization modes and enough
precision). These methods were called \emph{accurate} in
\cite{LlaveR90a}. To show that a method is accurate, one needs to
perform some mathematical analysis that ensures that the functions
involved have domains that are completely covered by the
discretization (e.g. if one discretizes by power series, one would
need to show that the functions involved have domains which are
circles).

In this paper, we have gone even further in parameters by relying on
the standard methods of numerical analysis (e.g. reruns, increasing
accuracy and checking it does not change).

\subsubsection{Some other remarks}
We conclude by mentioning that, in different problems,
\cite{FiguerasH12} developed a similar methodology to compute reliably
close to breakdown based on efficient algorithms backed up by
a-posteriori theorems. This was used in \cite{FiguerasH15} to discover
different scenarios for breakdown of normal-hyperbolicity.

As we will see, we have obtained what seems a new scenario of
breakdown.  Investigating the domain of universality will take us far
from our main goal and will be postponed.

\subsection{Organization of the paper}
The spin-orbit model is presented in Section~\ref{sec:model} and its
Poincar\'e map is described in Section~\ref{sec:map}. Our approach
works with any Diophantine frequency, but in the numerical experiments we have
considered two frequencies, namely the golden ratio and another
irrational number related to the golden mean but closer to $1$ --the
resonance--.  Both numbers are {\em noble} (the continued fraction
expansion is $1$ after a finite number of terms).

We have computed the so-called basins of rotation number, which give a
qualitative picture of the dynamics as well as information on the
value of the drift parameter, for given values of the dissipation (see
Section~\ref{sec:rotation}).

The computation of the tori and their exploration near breakdown is
presented in Section~\ref{sec:breakdown}, while
Section~\ref{sec:studybreakdown} presents the scaling theory
of breakdown of tori (a growing part of the theory has
been made rigorous by now).

The construction of the invariant bundles and an exploration of the
angle between the stable and tangent bundles, is given in
Section~\ref{sec:bundles}. Moreover, Section~\ref{sec:domains} shows
the width of analyticity domains of the tori and the angle between the
invariant bundles.

The computation of the breakdown threshold through Sobolev's criterion
is given in Section~\ref{sec:breaksob}.

The links with renormalization theory and the study of scale-invariant
properties of the tori are studied in
Section~\ref{sec:breakdownscaling}.

Finally, some conclusions are presented in
Section~\ref{sec:conclusions}.

\section{The spin-orbit problem with tidal effects}\label{sec:model}

In this section, we describe the physical motivation for the model we
will investigate; even if not strictly needed for the mathematical or
numerical study, the physical interpretation motivates our questions
and the ranges of the parameters involved. We will quickly review the
most important features, but leave a more detailed derivation to the
referenced literature (\cite{Peale,Celletti2010}).

The spin-orbit problem describes a simplified model for the rotational
dynamics of a satellite, say $\mathcal{S}$, orbiting around a central
planet, say $\mathcal{P}$, and rotating around an internal spin-axis
(\cite{Beletsky,Celletti90I,Celletti90II,LaskarC,Wisdom}).  We assume
that the satellite is a triaxial ellipsoid with principal moments of
inertia $\mathcal{A} < \mathcal{B} < \mathcal{C}$.  Moreover, we make
the following simplifying assumptions:

\begin{enumerate}
\renewcommand*{\theenumi}{\sl A\arabic{enumi}}
\renewcommand*{\labelenumi}{\theenumi.}
\renewcommand*{\itemsep}{.5em}
\item \label{assump1} The barycenter of the satellite $\mathcal{S}$
  moves on an elliptic Keplerian orbit with the planet $\mathcal{P}$
  in one focus; the Keplerian orbit is characterized by a semimajor
  axis $a$ and an eccentricity $\ecc$.

\item \label{assump2} The satellite spin-axis coincides with the
  direction of the smallest physical axis of the ellipsoid.

\item \label{assump3} The spin-axis is perpendicular to the orbital
  plane.

\item \label{assump4} The satellite is assumed to be non-rigid, hence
  subject to a tidal torque.

\item \label{assump5} The tidal torque is a linear function of the
  angular velocity of rotation (\cite{macdonald,Peale}, see also
  \cite{CellettiL2014,LaskarC}).
\end{enumerate}
The unit of time is chosen to normalize the orbital period $T_{orb}$
to $2\pi$, which implies that the mean motion $n=2\pi/T_{orb}$ is
equal to one.

The equation of motion of the spin-orbit problem is introduced as
follows. Let $x$ be the {\em rotation angle} formed by the largest
axis of the ellipsoid (belonging to the orbital plane, due to
\ref{assump2} and \ref{assump3}) and the periapsis line. Then, the
spin-orbit problem under the conditions \ref{assump1}-\ref{assump4} is
described by the equation
\begin{equation}\label{eq.spin-diss}
\frac{d^2x}{dt^2}(t) + \eps \left(\frac{a}{r(t)}\right)^3 \sin \bigl(2
x(t) - 2f(t)\bigr) = \mathcal{T} _d \biggl(\frac{dx}{dt}(t),t
\biggr)\ ,
\end{equation}
where the {\em perturbative parameter} $\eps \bydef \frac{3}{2}
\frac{\mathcal{B}-\mathcal{A}}{\mathcal{C}}$ measures the equatorial
ellipticity of the satellite. The terms $r(t)=r(t;\ecc)$ and
$f(t)=f(t;\ecc)$ are known functions depending on the eccentricity and
related to time by means of the {\em Kepler equation}. In fact,
denoting by $t _0$ the initial time, Kepler equation $nt+t
_0=u-\ecc\sin u$ gives the eccentric anomaly $u$ as a function of the
time $t$, while $r$ and $f$ are related to $u$ by
\begin{equation} \label{eq.rtanf}
  r=a(1-\ecc \cos u)\ ,\qquad
  \tan \frac{f}{2} =\sqrt{\frac{1+\ecc}{1-\ecc}} \tan \frac{u}{2}\ .
\end{equation}

The right hand side of \equ{eq.spin-diss} represents the dissipative
effect, which is given by the tidal torque $\mathcal{T} _d=\mathcal{T}
_d(\frac{dx}{dt}(t),t)$; using the model developed in
\cite{macdonald,Peale} (see Assumption~\ref{assump5}), $\mathcal{T}
_d$ takes the form
\begin{equation}
\label{eq.mcdonald}
\mathcal{T} _d \biggl(\frac{dx}{dt}(t), t\biggr) \bydef - \dis
\left(\frac{a}{r(t)}\right)^6 \biggl(\frac{dx}{dt}(t) - \frac{d
f}{dt}(t) \biggr)\ ,
\end{equation}
where $\dis > 0$ is called the {\em dissipative constant}, which
depends on the physical features of the satellite (density, rigidity,
etc.).

We remark that the astronomical values of $\eps$ are small for many
satellites and planets of the Solar system, typically of the order of
$10^{-4}$ for celestial bodies such as the Moon and Mercury, and that the
dissipative term is typically small with respect to the conservative
term, say $O(\eps^2)$.

For $\dis=0$, equation \equ{eq.spin-diss} is conservative, and for
$\eps=\dis=0$ it is integrable. Thus, for small values of $\eps$ and
$\dis$, we obtain a nearly-integrable and nearly-Hamiltonian system.

\vskip.1in

Taking the average of the tidal torque over one orbital period (see,
e.g., \cite{Peale,LaskarC}), equation \equ{eq.spin-diss} becomes
\begin{equation}
\label{eq.avg-spinxy}
\frac{d^2x}{dt^2}(t) + \eps \left(\frac{a}{r(t)}\right)^3 \sin
\bigl(2x(t)-2f(t) \bigr) = -\dis \biggl(\bar
L(\ecc)\frac{dx}{dt}(t)-\bar N(\ecc)\biggr)\ ,
\end{equation}
where
\begin{equation} \label{LNave}\begin{split}
    \bar L(\ecc)&= {\frac{1}{(1-\ecc^2)^{9/2}}}
    \biggl(1+3\ecc^2+\frac{3}{8}\ecc^4 \biggr) \ , \\
    \bar N(\ecc)&= \frac{1}{(1-\ecc^2)^6}
    \biggl(1+\frac{15}{2}\ecc^2+\frac{45}{8}\ecc^4+ \frac{5}{16}\ecc^6
    \biggr)\ .
\end{split} \end{equation}

The averaged model is obtained from the non-averaged model by
eliminating the high-frequency terms in the tidal torque; notice that
the terms in \eqref{LNave} are exact and they are not truncation of
expansions in the eccentricity.

If we write the second order equation \eqref{eq.spin-diss} as a first
order system in phase space, we obtain:
\begin{equation} \label{eq.spinxy}
\begin{split}
  \dot x&= y\\
  \dot y&=- \eps \left(\frac{a}{r(t)}\right)^3 \sin\bigl(2
  x-2f(t)\bigr)-\dis \left(\frac{a}{r(t)}\right)^6 (y - \dot
  f(t)) \\
  \dot t&=1
\end{split}
\end{equation}
and similarly for equation \equ{eq.avg-spinxy}:
\begin{equation} \label{A1}
  \begin{split}
    \dot x&= y \\
    \dot y&=- \eps \left(\frac{a}{r(t)}\right)^3 \sin(2 x-2f(t))-\dis
    \biggl(\bar L(\ecc)y-\bar N(\ecc)\biggr) \\
    \dot t&=1\ .
\end{split}
\end{equation}
Equations \eqref{eq.spin-diss} and \eqref{eq.avg-spinxy} are defined
over the phase space $[0,2\pi)\times \mathbb{R}$, which can be endowed
  with the standard scalar product and the symplectic form $\Omega =
  dy \wedge dx$.
Both equations \eqref{eq.spin-diss} and \eqref{eq.avg-spinxy} are {\em
  dissipative} in the sense that the phase space volume contracts
under the evolution of the flow. Hence, if $J_t$ is the determinant of
the linearized flow, Abel's formula gives
\begin{equation*}
  \det(J_t) = \exp{\int_0^t \Tr(A(s))\ ds}\ ,
\end{equation*}
where $\Tr$ denotes the trace and $A(t)$ is the differential of the
vector field at the time $t$ flow.
% $\bydef D{\underline
%   f}({\underline F}_t(y,x,t))$ with $\underline f$ denoting the vector
% field and ${\underline F}_t={\underline F}_t(y,x,t)$ the flow.

With reference to equation \equ{eq.spinxy}, we obtain
\begin{equation*}
  \Tr(A(t))=-\dis \left(\frac{a}{r(t)}\right)^6
\end{equation*}
and consequently the volume contracts as
\begin{equation}
  \label{eq.Jtxy}
  \det(J_t)=\exp\biggl(-\dis \int_0^t \left(\frac{a}{r(t)}\right) ^6\,
  ds \biggr)\ .
\end{equation}
For the averaged tidal torque model \equ{A1}, we obtain
$$
\Tr(A(t))=-\dis \bar L(\ecc)
$$
and consequently the volume contracts as
\begin{equation}
  \label{eq.Jtavgxy}
  \det(J_t) = \exp\bigl(-\dis \bar L(\ecc)\ t\bigr)\ .
\end{equation}

It is not difficult to see by an explicit computation (compare with
\cite{CCGL20a}), that at $t = 2\pi$ both \eqref{eq.Jtxy} and
\eqref{eq.Jtavgxy} give the same asymptotic contraction rate of the
volume.

Note that, when we do continuation in $\eps$, the value of the
dissipation changes. This is somewhat different from previous studies
in which the continuation was done for families of constant
dissipation.

\begin{rem}\label{averages}
The relation between the averaged and non-averaged models has been
considered puzzling in the literature. The difficulty is that
averaging methods are usually justified for short time.  On the other
hand, attractors involve the evolution over very long times.

As we will see in Section~\ref{sec:map}, one of the advantages of the
approach of \cite{CCGL20a} is that the computations of the attractors
are based on time-1 maps (for which the standard justification of the
averaging applies). The a-posteriori format of the result in
\cite{CallejaCdlL2013} shows that, when the maps are close,  the attractors of the maps are
close.  Hence, the a-posteriori results, that
justify an approximations for finite time, justify results for objects over all times.

Besides the general argument above, \cite{CCGL20a} also presents a
self-contained elementary argument showing that attractors of the two
models are close over all times.  In this paper, we present some
numerical explorations that show that the rigorous results on the
similarity of the attractors of the two models can be observed in numerical experiments, see
Figure~\ref{fig.rotx} and Figure~\ref{fig.avg-navg}.
%%%%%% The simularity remains also in much more delicate figures and it is, therefore omitted.
\end{rem}

\section{The Poincar\'e spin-orbit map}
\label{sec:map}

To perform the computations that will be presented in the next
sections, it is convenient to reduce the study of \equ{eq.spin-diss}
to the discrete system associated to the vector field
\eqref{eq.spinxy} by a return map.
According to the procedure detailed in \cite{CCGL20a}, we consider the
Poincar\'e map $P _\ecc$ associated to \equ{eq.spinxy}, which is
obtained as follows. We introduce the map $P_ \ecc$ as
\begin{equation}
  \label{eq.Pe}
  P _\ecc(x_0,y_0;\eps) \bydef
  \begin{pmatrix}
    x(2\pi;x_0,y_0,\eps) \\ y(2\pi;x_0,y_0,\eps)
  \end{pmatrix}\ ,
\end{equation}
where $x(2\pi;x_0,y_0,\eps) $ and $y(2\pi;x_0,y_0,\eps)$ represent the
solution of \eqref{eq.spinxy} at time $t=2\pi$ with initial conditions
$(x_0,y_0)$ at $t=0$. Writing $P _\ecc$ in components as $P _\ecc
\equiv (P _\ecc^{(1)},P _\ecc^{(2)})$, the \emph{spin-orbit Poincar\'e
map} is given by
\begin{equation}
  \label{SOmap}
  \begin{split}
  \bar x&= P _\ecc^{(1)}(x,y;\eps)\ ,\\ \bar y&= P
  _\ecc^{(2)}(x,y;\eps)\ .
 \end{split}
\end{equation}
For computational reasons, it is convenient to consider the change of
coordinates
\begin{equation}
  \label{eq.xy2bg}
  \Psi _\ecc \bydef 2\pi
  \begin{pmatrix}
    1 & 0 \\ 0 & 1 - \ecc
  \end{pmatrix}
\end{equation}
and define the map
\begin{equation} \label{map2}
  G_\ecc \bydef \Psi _\ecc \circ P _\ecc \circ \Psi_\ecc^{-1}\ ,
\end{equation}
which is in fact the Poincar\'e spin-orbit map for \eqref{eq.spinxy}
with the change of coordinates given by Kepler's equation, $t = u -
\ecc \sin u$.

The map \equ{SOmap}, equivalently \equ{map2}, inherits some properties
from the differential equation \equ{eq.spin-diss} (equivalently,
\eqref{eq.spinxy}). Notably, \equ{SOmap} (or \equ{map2}) is a
nearly-integrable dissipative map, which is indeed {\em conformally
  symplectic} according to the following definition (see
\cite{CallejaCdlL2013}).

\begin{defn} \label{def:conformallysymplectic}
Let $\mathcal{M} = \T^n\times U$ with $U\subseteq \mathbb{R}^n$ an
open, simply connected domain with smooth boundary and with symplectic
form $\Omega$. A family of diffeomorphisms $f _\vartheta\colon
\mathcal M\rightarrow\mathcal M$, depending on a parameter
$\vartheta$, is conformally symplectic, if there exists a function
$\lambda \colon \mathcal{M}\to\mathbb{R}$ such that
\begin{equation}\label{CS}
  f _\vartheta^* \Omega = \lambda\Omega\ ,
\end{equation}
where $f _\vartheta^*$ denotes the pull--back.
\end{defn}

The quantity $\lambda$ is called the conformal factor. For $\lambda=1$
we recover the symplectic case. For $n=1$ any diffeomorphism is
conformally symplectic with $\lambda(x)=\pm|\det(D f
_\vartheta(x))|$. For $n\ge 2$ it is known that $\lambda$ is constant
(see \cite{CallejaCdlL2013} for details).

For later use, we denote by $J$ the symplectic matrix representing
$\Omega$ at ${\underline z}$, namely $\Omega_{\underline
  z}({\underline u},{\underline v})=({\underline u},J({\underline
  z}){\underline v})$ for ${\underline u},{\underline
  v}\in\mathbb{R}^n$ and $(\cdot,\cdot)$ is the standard inner product
in $\R^n$.  For the spin-orbit map, the matrix $J$ is given by
\begin{equation}
 \label{matrixJ}
 J =
 \begin{pmatrix}
  0 & 1 \\ -1 & 0
 \end{pmatrix}\ .
\end{equation}

\vskip.1in

The Poincar\'e spin-orbit map \equ{SOmap} (or \equ{map2}) is
conformally symplectic (\cite{CCGL20a}) with the conformal factor
given by
\begin{equation} \label{lambda}
  \lambda = \exp \biggr(-\dis \pi \frac{3 \ecc^4+24 \ecc^2+8}{4
    \left(1-\ecc^2\right)^{9/2}}\biggl)\ .
\end{equation}
The system contracts the volume for $\dis>0$, expands the volume for
$\dis<0$, and it is neutral for $\dis=0$. In the following sections,
we will only take $\dis> 0$.

Note that the conformal factor of the averaged model
\eqref{eq.avg-spinxy} will have the same $\lambda $ given in
\eqref{lambda}, since at time $t=2\pi$ both contraction rates
\eqref{eq.Jtxy} and \eqref{eq.Jtavgxy} coincide precisely due to the $\bar L(\ecc)$ choice.
For the calculation leading to \eqref{LNave} we refer to
\cite[p. 15]{CCGL20a}.

\section{Basins of rotation numbers}
\label{sec:rotation}

In this section we present exploratory results based on the
computation of the rotation number of many orbits of the spin-orbit
problem. These computations provide a procedure to obtain an
approximate value of the drift parameter associated to a given
frequency of an invariant torus.  Such information is essential to
start the procedure detailed in \cite{CCGL20a} for the computation of
a KAM torus with fixed frequency.

The phase space associated to the spin-orbit map \equ{SOmap} is a
subset of the cylinder $[0,2\pi) \times \mathbb{R}$.  For numerical
  reasons, it is convenient to take the eccentric anomaly $u$ as
  independent variable, instead of time. This can be achieved by
  taking into account Kepler's equation, which provides a relation
  between $t$ and $u$: $ t = u - \ecc \sin u$.  More precisely, let
  $s(x;u,\ecc) \bydef \sin(2x(t)-2f(t))$, which depends on $u$ and
  $\ecc$ through $f$, that satisfies the following identities:
  \beq{ffue} \cos f = \frac{\cos u - \ecc}{1 - \ecc \cos u} \qquad
  \text{and} \qquad \sin f = \frac{\sqrt{1 - \ecc^2}\sin u }{1 - \ecc
    \cos u}\ .  \eeq We use simple trigonometric identities to write
  the function $s$ as
\begin{equation}
  \label{eq.sx}
  s(x; u, \ecc) = \sin(2x) (2\cos^2 f-1) - \cos(2x)2\cos f \sin f
\end{equation}
and we define the function $c$ as
\begin{equation}
  \label{eq.cx}
  c(x; u, \ecc) \bydef \cos(2x) (2\cos^2 f-1) + \sin(2x)2\cos f \sin
  f\ .
\end{equation}
Note that both quantities are easier to compute than the one in terms
of $\tan(f/2)$ in \eqref{eq.rtanf} which has singularities depending
on the value of $u$.  Moreover, \eqref{eq.sx} and \eqref{eq.cx}
satisfy the following relations, which will be used during the
integration with Taylor's numerical integration method:
\begin{equation*}
 \frac{\partial s}{\partial x}(x; u, \ecc) = 2c(x; u, \ecc)\ , \qquad
 \frac{\partial c}{\partial x}(x; u, \ecc) = -2 s (x; u, \ecc) \ .
\end{equation*}
As a consequence of the change of coordinates, we deduce that
\begin{equation*}
 \frac{d f}{dt} = \left(\frac{a}{r(u)}\right)^2 \sqrt{1 - \ecc ^2}\ .
\end{equation*}
The spin-orbit problem in \eqref{eq.spinxy} can then be expressed in
terms of the independent variable $u$ as
\begin{equation}
  \label{eq.spinbg}
  \frac{d^2 \beta}{du ^2}(u) - \frac{d\beta}{du}(u) \frac{a}{r(u)}
  \ecc \sin u + \eps \frac{a}{r(u)} s(\beta;u, \ecc) = - \dis
  \left(\frac{a}{r(u)}\right)^5 \biggl(\frac{d\beta}{du}(u) -
  \frac{a}{r(u)} \sqrt{1 - \ecc ^2}\biggr)\ ,
\end{equation}
where $\beta$ and $\gamma$ are defined as
\begin{equation*}
  \label{eq.x2b}
  \beta (u) \bydef x (u - \ecc \sin u)\ , \qquad \gamma (u) \bydef
  \frac{d\beta}{d u}(u) = \frac{r(u)}{a}y(u - \ecc \sin u)\ .
\end{equation*}
The time-($2\pi$) map associated to \equ{eq.spinbg} is related, up to
a factor, to the map defined in \eqref{map2}.

\subsection{Computation of the rotation numbers}

In this section, we present the computation of the rotation numbers
obtained by a direct iteration; in Section~\ref{sec:basins}, we will
map the regions covered by different attractors. First, we recall that
the definition of the rotation number associated to equation
\equ{eq.spinxy} (equivalently \equ{A1}) is:
\begin{equation}
  \label{eq.rotnum}
  \rho = \lim _{n \to \infty} \frac{1}{n} \sum _{j = 1}^{n-1} y(2\pi
  j)\ .
\end{equation}
The rotation number of orbits of two-dimensional maps could fail to
exist, but when the orbits lie in an invariant circle (the main object
of our interest in this paper) such rotation number exists.  Also, all
the orbits in the basin of attraction of a circle have the same
rotation number.

The convergence of the limit in \equ{eq.rotnum} may be very slow,
say of the order $O(1/n)$.  As shown in \cite{Sander2017}, a
method alternative to \equ{eq.rotnum} to compute the rotation
number is through the following formula:
\begin{equation}
  \label{rho}
  \rho = \lim _{n \to \infty} \biggl[\frac{1}{\sum _{j = 1}^{n-1}
      \phi(\tfrac{j}{n})} \sum _{j = 1}^{n-1} y (2\pi j) \phi
    (\tfrac{j}{n})\biggr]\ ,
\end{equation}
where $\phi$ denotes a weight function, e.g.,
\begin{equation}
  \label{eq.weightfun}
  \phi(z) \bydef \exp \biggl(\frac{-1}{z^2(1-z)^2} \biggr), \qquad z
  \in (0,1)\ .
\end{equation}
It is shown in \cite{Sander2017} that, if the rotation number is
sufficiently irrational (and the corresponding circle is smooth), the
limit in \eqref{rho} is reached very fast. For example, if the circle
is analytic and its frequency is Diophantine, the limit is reached
super-exponentially fast.

In contrast, if the limiting rotation number is rational (or if the
circle is not smooth), the convergence of \eqref{rho} is slow.  Hence,
the speed of convergence on \equ{rho} can be used as a diagnostic
whether the torus is present and smooth.  For our purposes, the main use of
\eqref{rho} is to adjust parameters, so that an attractor with the
right rotation number is found.

The computation of the rotation number for the system
\eqref{eq.spinxy}, associated to an attractor close to the initial
condition $(x _0, y _0) \in [0,2\pi) \times \mathbb{R}$, can be
  performed according to Algorithm~\ref{alg.rotx}, described
  below\footnote{In Algorithm~\ref{alg.rotx}, the symbol $\gets$ means
  assignment, that is, the quantity on the right gives the quantity on
  the left.}.

\begin{alg}
\label{alg.rotx}
 Computation of the rotation number for the spin-orbit problem with
 tidal torque, see equation \equ{eq.spinxy}.
 \begin{enumerate}
\setlength{\itemsep}{.8em}
\renewcommand*{\theenumi}{\emph{\arabic{enumi}}}
\renewcommand*{\labelenumi}{\theenumi.}
  \item [$\star$] {\tt Inputs:}
  \begin{description}
   \item [Initial condition] $(x _0, y _0) \in [0,2\pi)\times
     \mathbb{R}$.
    \item [Spin-orbit parameters] $\eps > 0$, $\ecc \in [0,1)$, $\dis
      > 0$, and $\lambda$ as in \eqref{lambda}.
    \item [Positive integers] $\delta$, $n_1$, $n_2$ with $n _1 < n
      _2$, and $n _0 \gets \lceil -14/\log_{10} \lambda \rceil$.
    \item [Weight function] $\phi$, for instance, the one given in
      \eqref{eq.weightfun}.
  \end{description}
  \item [$\star$] {\tt Outputs:} Approximate rotation number $\rho$
    for \eqref{eq.spinxy} (equivalently, for \equ{A1}).
  \item $(\beta _0, \gamma _0) \gets (x _0,y _0 / (1-\ecc))$.
  \item \label{alg.rotx.trans} $(\beta _0, \gamma _0) \gets G_\ecc ^{n
    _0}(\beta _0, \gamma _0)$ by integrating \eqref{eq.spinbg} for
    $2\pi n _0$ times in $[0,\pi)\times \mathbb{R}$.
  \item Store $\gamma _j \gets \gamma(2\pi j; \beta _0, \gamma _0)$
    for $j = 1, \dotsc, n _2$ by integrating \eqref{eq.spinbg}.
  \item \label{alg.rotx.loop} $n \gets n _1$.
  \item \label{alg.rotx.loop2} $ \rho \gets \bigl(\sum _{j = 1}^{n-1}
    \phi (\tfrac{j}{n})\bigr)^{-1}\: \sum _{j = 1}^{n-1} \gamma _j \phi
    (\tfrac{j}{n})$.
  \item Iterate steps \eqref{alg.rotx.loop}-\eqref{alg.rotx.loop2}
    with $n \gets n + \delta$ until $\rho$ does not vary or $n > n
    _2$.
  \item If convergence, return $\rho/ (1 - \ecc)$.
 \end{enumerate}
\end{alg}

To find a reliable rotation number of a randomly chosen orbit, it is
convenient to discard a transient number of iterations to ensure that
the orbit has converged to the 1-D attractor; the parameter $n_0$ in
step \ref{alg.rotx.trans} of Algorithm~\ref{alg.rotx} implements
this. Unfortunately, the time that an orbit needs to settle in an
attractor is hard to estimate and it may be surprisingly large (see
\cite{CellettiL2014}), because an orbit may be entertained near an
attractor for a long time before landing in the final one.  This
effect is very hard to predict and changes wildly with the initial
conditions.

Since the only harm in overestimating the transition is the use of
more computer time, we have tried to overestimate it. Note also that
iterating different points is very paralellizable.

\begin{figure}[ht]
\resizebox{.95\columnwidth}{!}{\input{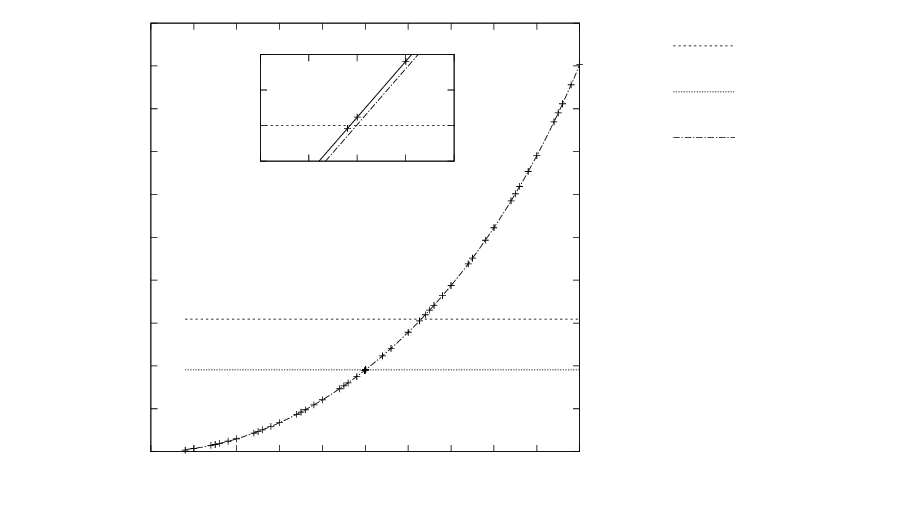}}
 \caption{Rotation numbers ({\tt +}) associated to \equ{eq.spinxy}
   from Algorithm~\ref{alg.rotx} in terms of the eccentricity $\ecc$
   with weight function as in \eqref{eq.weightfun}, $n _1 = 4500$, $n
   _2 = 4600$, $\delta = 10$, $\eps=10^{-4}$, $\mu=10^{-5}$ and $(x
   _0, y _0) = (0, 1.38)$.  The fixed frequencies $\omg _1$ and $\omg
   _2$ are given in \eqref{eq.omg1} and \eqref{eq.omg2}, the predicted
   rotation number for the averaged spin-orbit is given by $\bar
   N(\ecc)/\bar L(\ecc)$.  Finally a zoom-in is shown nearby
   $\omega_2$ with the {\tt +} connected by a straight line.}
 \label{fig.rotx}
\end{figure}

We report in Figure~\ref{fig.rotx} the values of the eccentricities
that lead to invariant attractors with given rotation number. As
sample cases, we are interested to the following two irrational
numbers:
\begin{align}
 \omg _1 &\bydef \gamma _g^+ \approx
 \mathtt{1.618033988749894848}\ldots , \label{eq.omg1} \intertext{and}
 \omg _2 &\bydef 1 + \frac{1}{2+\gamma _g^-} \approx
 \mathtt{1.3819660112501051} \ldots \label{eq.omg2}
\end{align}
with $\gamma _g^\pm$ equal to the golden ratio and its conjugate:
$\gamma _g^\pm \bydef \frac{\sqrt{5}\pm 1}{2}$.  Note that both
$\omg_1, \omg_2$ are \emph{noble} numbers (i.e.  the continued
fraction expansion has only $1$ after a certain point).  It is known
that one of the predictions of the standard renormalization group
\cite{MacKay82} is that for maps in the universality domain many of the fine properties at breakdown are
similar for all tori with noble rotation numbers.

Figure~\ref{fig.rotx} shows $\omg _1$ and $\omg _2$ with two different
dashed lines. Moreover, it shows with a dotted-dashed line the term
$\bar N(\ecc)/\bar L(\ecc)$ (see \eqref{LNave}), which is the
predicted rotation number in the model with the averaged tidal torque,
see \cite{ARMA, CCGL20a}. Extremely close to the dotted-dashed line,
the crosses are obtained by computing the rotation number for the
non-averaged model \equ{eq.spinxy}. To highlight the proximity between
the results associated to the averaged and the non-averaged models,
Figure~\ref{fig.rotx} shows also a zoom-in nearby the frequency $\omg
_2$.  This shows that the rigorous results in \cite{CCGL20a}
justifying the averaging method give very good numerical values.

\subsection{Computation of the basins of rotation numbers}\label{sec:basins}

A qualitative picture of the dynamics can conveniently be achieved by
computing the {\em basins of rotation numbers}. These basins are
obtained by applying recurrently Algorithm~\ref{alg.rotx}; given a
grid in $(x _0, y _0)$ of size $n _x \times n _y$ (for some positive
integers $n_x$, $n_y$) in a window of the cylinder $[0,2\pi) \times
  \mathbb{R}$, we apply the Algorithm~\ref{alg.rotx} and represent the
  rotation number as a color map, using interpolation for the elements
  that are not in the finite grid.

Notice that the iterations of each orbit are completely independent
from each other, which makes the procedure completely
parallelizable. Each orbit could be assigned to a thread (using
e.g. {\sc openMP} \cite{Dagum1998openmp}) or, working on different
regions, could be distributed using coarse grained distributors such
as GNU-parallel or HTCondor (see \cite{Tange2011a,htcondor}). We have
used these techniques systematically.

Although we have not yet used GPU's, the computation of these basins
are indeed suitable for GPU's. For preliminary calculations, the
single precision available even in consumer grade GPU's is perfectly
acceptable. However, when dealing with the more detailed calculations
of tori, even double precision is not enough and we need extended
precision.

\begin{figure}[ht]
\resizebox{.98\columnwidth}{!}{\input{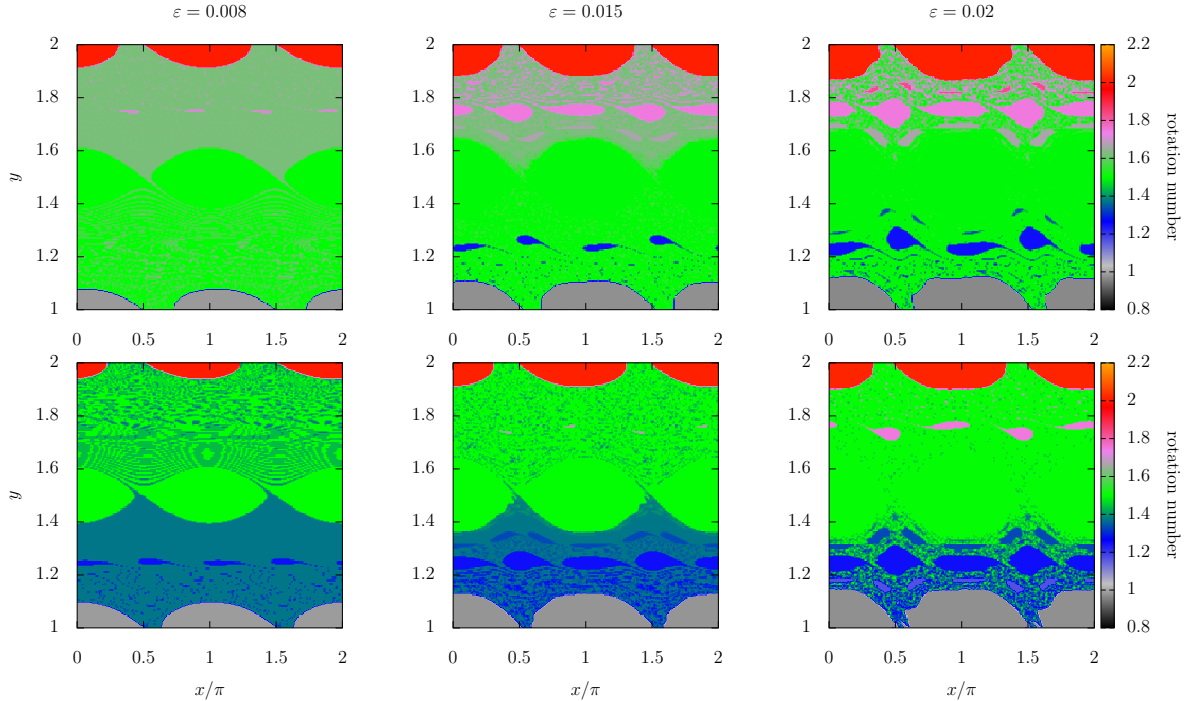}}
 \caption{Basins of rotation numbers with grid size $500\times 500$
   and window $[0,2\pi) \times [1,2]$ for the spin-orbit problem
     \eqref{eq.spinxy} with $\dis=10^{-3}$, different values of $\eps$
     per each column, and $\ecc = 0.316$ in the upper row,
     $\ecc=0.247$ in the lower row.  The color scale provides the
     value of the rotation number for each initial condition $(x,y)$,
     using interpolation for points outside the grid.}
 \label{fig.basins}
\end{figure}

Figure~\ref{fig.basins} shows some basins of rotation numbers for
different values of the parameters. On average, each of the plots,
using the {\tt taylor} program \cite{JorbaZ2005}, needs around 12
minutes using 143 CPUs.  The color scale gives the rotation number for
every initial condition.  We provide the results for two different
values of the eccentricity.  The plots clarify regions of librational
and rotational motion. At a finer inspection, the regions close to the
location of the invariant tori with frequencies $\omega _1$ in the
upper plots of Figure~\ref{fig.basins} and $\omega _2$ in the lower
plots, become more irregular as the perturbing parameter $\eps$ ranges
over values far below from what will be the critical threshold (left
column), close to the critical threshold (middle column) or far above
the critical threshold (right column).  This behavior is consistent
with the theory on the existence of KAM tori developed in
\cite{CallejaCdlL2013}.

\section{Explorations of the breakdown}\label{sec:breakdown}

A rotational torus defined for a conformally symplectic system (see
Definition~\ref{def:conformallysymplectic}) is introduced as follows.

\begin{defn}\label{defn.confotori}
 Let $f _\vartheta \colon \mathcal{M} \to \mathcal{M}$ be a
 conformally symplectic system and let $\omg$ be an irrational
 number. A map $K \colon \T^n \to \mathcal{M}$ is said to be a
 rotational torus for $f _\vartheta$ with frequency $\omg$, when there
 exists $\vartheta _0$ such that
 \[
  f _{\vartheta _0} \circ K(\theta) = K (\theta + \omg) \qquad
  \text{for all } \theta \in \T^n.
 \]
\end{defn}

The algorithm in \cite{CallejaCdlL2013} deals with tori whose
frequency has Diophantine properties and it computes iteratively
corrections for the map $K$ and the parameter $\vartheta _0$ such that
Definition~\ref{defn.confotori} is satisfied up to a numerical
tolerance. If there is an extra parameter in the problem, it is
standard to use the Newton method to implement a continuation
method. The solutions of the invariance equation for a parameter value
are taken as an initial guess for the problem with an slightly larger
parameter value and applying the Newton method, we obtain a solution
for the increased parameter value.

In many situations, this continuation method stops because the
algorithm becomes unreliable (e.g the expansions used in the
truncation of the problem do not reach the function) or because the
torus ceases to exist (as a smooth object).

In the following sections we present very accurate computations of
these tori and we discuss some interesting phenomena near the
spin-orbit breakdown.

\subsection{Computation of the tori and their continuation towards the breakdown}

The theorems and algorithms in \cite{CallejaCdlL2013} provide a way to
efficiently compute KAM tori in conformally symplectic systems. Those
results have been adapted for the spin-orbit problem in
\cite{CCGL20a}. The main difference between the symplectic and the
conformally symplectic cases is that, in the dissipative case at each time that the torus is
corrected for a fixed rigid rotation, a parameter of the system must
be corrected as well. In fact, while the torus with a fixed frequency
can exist for each parameter value in a conservative system, in a dissipative one the torus may only exist for certain
parameter values. Moreover, in a conformally symplectic system a torus
with a given rotational frequency exists for parameter values belonging to a
whole interval, the so-called Arnold tongue (\cite{Arnold65,CCFL14}).

The method detailed in \cite{CallejaCdlL2013} is really efficient and
accurate in computing the torus and the drift parameter of the system
given approximate values. This makes it extremely suitable for a
continuation algorithm (we take the solution for a parameter value as
an approximate solution for a slightly bigger value of the
parameter\footnote{As a practical matter, we use some extrapolation
from previous values.}).

Furthermore the continuability argument in \cite{CallejaCdlL2013}
shows that the continuation method will get as close as desired to the
boundary of existence if given enough computer resources (getting
extremely close to the boundary of existence may require using a very
large number of Fourier modes, extended precision arithmetic,
etc.). Note that several KAM proofs and their attendant algorithms do
not satisfy this property. This requires that the parametrization used
reaches the maximal domain, see \cite{LlaveR90a}).

In reality, of course, one has limited computer resources (limited
memory, computer time, and time of the programmer).  In what follows,
we will explore the limits obtained using standard computers (desktops
or small servers) at the time of the writing.  We note that the
algorithms based on automatic reducibility specified in
\cite{CallejaCdlL2013} are very efficient both in memory required, in
number of operations and in the dimension of the objects studied. The
algorithm is also based on very structured operations which make it
easy to change the precision of the arithmetic and the number of
Fourier modes in a modern programming language.

One of the results of our explorations is that the main bottleneck for
the accuracy of the calculation near the breakdown is the precision of
the arithmetic. Therefore, in our studies, we use high precision
arithmetic and a very high order ODE integrator.  Both of these
improvements are doable in a reasonable programming time, thanks to
the {\tt taylor} program \cite{JorbaZ2005}, which generates with
minimal programming effort a Taylor solver (very high order) with
extended precision arithmetic. This is crucial to obtain high
precision for the Poincar\'e spin-orbit map.

\subsection{Some implementation details; running times}
For the frequency \eqref{eq.omg1} the computation was done in a
machine with Intel Xeon Gold 5220 CPU at 2.20GHz with 18 CPUs with
hyperthreading which simulates 36 CPUs, while for \eqref{eq.omg2} a
machine with Intel Xeon Gold 6154 CPU at 3.00GHz with 74 CPUs with
hyperthreading which simulates 148 CPUs. In case of parallelism, we
have always requested 32 CPUs.

\subsubsection{Choice of continuation parameters}
The procedure in \cite{CCGL20a} involves computing the derivatives of
the spin-orbit Poincar\'e map \eqref{SOmap} with respect to $(x,y)$
and with respect to the parameter to be corrected. In our case, we
will use the eccentricity as the parameter to be adapted. The
eccentricity has a clear physical interpretation. We will perform our
calculations for tori of the two frequencies \eqref{eq.omg1} and
\eqref{eq.omg2}. Thus, we get $(K _{\omg _1}, \ecc _{\omg _1})$ and
$(K _{\omg _2}, \ecc _{\omg _2})$ such that
\begin{equation}\label{eq.inveq}
 P _{\ecc _{\omg _j}} \circ K _{\omg _j} (\theta) = K _{\omg _j}
 (\theta + \omg _j)\ , \qquad j \in \{1,2\} \ .
\end{equation}
We fix the dissipation parameter to $\dis =10^{-3}$ and we perform a
continuation w.r.t. the parameter $\eps$.  The continuation is based
on using a cubic extrapolation from the previous results to provide an
initial guess for the next continuation step and then, polishing this
guess by running the Newton method (as many times as needed to achieve
the desired accuracy, see below).

\subsubsection{Step control and adaptability of the algorithm}

Besides the usual step control in continuation methods, we incorporate
some control flags to ensure that the computed embedded torus
satisfies the invariance equation with enough accuracy and it is a
reasonable function.

One of the flags monitors the error in the invariance equation
\eqref{eq.inveq}. To declare the computation successful, we require
that the error is below a certain value.

We note that in this problem, specially near the breakdown, it is
possible to obtain spurious solutions that indeed solve very
accurately the invariance equation, but which are very unreasonable
functions. Hence, we introduce a second flag which checks that the
tail of Fourier coefficients has norms between two tolerances
$\epsilon _1 \leq \epsilon _2$. If the size of the tail is bigger than
$\epsilon_2$, we increase the number of coefficients used; if the tail
is smaller than $\epsilon _1$, we decrease the number of Fourier
coefficients (to increase the speed). Hence, we declare a computation
successful when it achieves small residual and when the solution has a
small enough tail, so that it can be computed with different levels of
truncation.

In case of failure, in any of the control flags, remesh and increase the number of the Fourier coefficients in the last
successful computation, run the Newton method to get more accuracy
(and check that the new solution has good flags) and repeat the
process. Of course, there is a global limit to the number of Fourier
modes and the calculation finishes when we go over the limit. In this
run, we have used $L _{\max} = 65536$.
Figure~\ref{fig.stats-avg-navg-mesh} shows the different mesh values
used during the continuation. We also use standard adaptative step
control in continuation, which stops the calculation when a step below
the minimum does not achieve that the Newton method converges. If the
Newton method succeeds, the initial guess for the next step will be
obtained by doing a cubic extrapolation of the previously computed
tori with the same meshsize.

\subsubsection{Run times}
We have made several runs with different values for these tolerances
and checked that the results do not change appreciably by changing
the choices of parameters of the algorithm.
\begin{figure}[ht]
% \vglue2cm
% \hglue-3cm
\resizebox{\columnwidth}{!}{\input{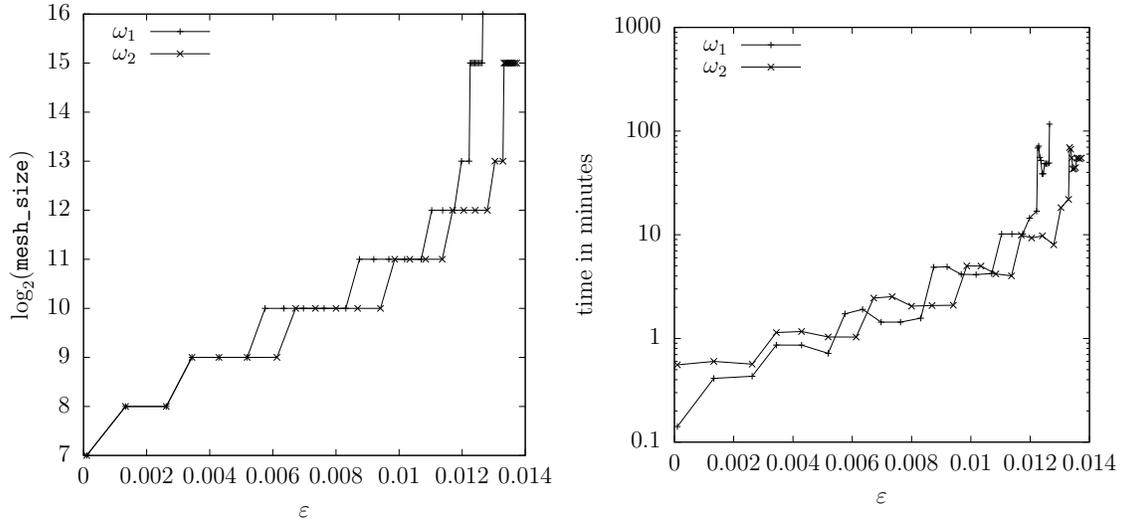}}
 \caption{Mesh size and execution time (in log10 minutes) of
   successful continuation steps for the non-averaged spin-orbit
   problem \eqref{eq.spinxy} and for the frequencies \eqref{eq.omg1}
   and \eqref{eq.omg2}.}
 \label{fig.stats-avg-navg-mesh}
\end{figure}

For the spin-orbit model, we use $70$ digits of precision, $10^{-35}$
of Newton tolerance for the error of the invariance equation
\equ{eq.inveq}, tail tolerance between $\epsilon _1=10^{-55}$ and
$\epsilon _2=10^{-28}$, and $10^{-70}$ as tolerance for the absolute
and relative errors in the Taylor integration.

\begin{figure}[ht]
\resizebox{\columnwidth}{!}{\input{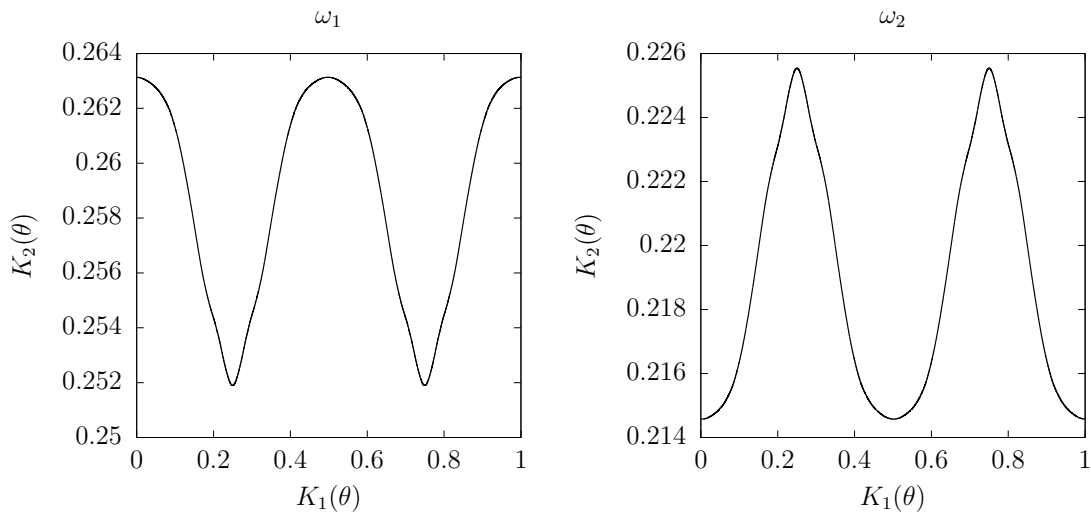}}
 \caption{Tori at the numerical breakdown for the non-averaged \eqref{eq.spinxy} spin-orbit problems and for the frequencies
   \eqref{eq.omg1} (left panel) and \eqref{eq.omg2} respectively.}
 \label{fig.criticaltori}
\end{figure}

The last values for the perturbing parameter $\eps$, that we reached
using the $\eps$ values in the continuation strategy described above
are \eqref{eq.epsomg1} for $\omg _1$ and \eqref{eq.epsomg2} for $\omg
_2$:
\begin{align}
 \eps _{\omg _1} ^{c} &=
 \texttt{1.265364670507455833901687945901589e-02}
 \ , \label{eq.epsomg1}\\ \eps _{\omg _2} ^{c} &=
 \texttt{1.372784208166277584502189850802202e-02}
 \ . \label{eq.epsomg2}
\end{align}
The corresponding tori are plotted in Figure~\ref{fig.criticaltori} and the values of the drift parameter for $\omg _1$ and
$\omg _2$ are:
\begin{align}
 e_{\omg _1} ^{c} &=
 \texttt{3.1701530650181080344148405746028386e-01}
 \ , \label{eq.eomg1}\\
 e_{\omg _2} ^{c} &=
 \texttt{2.4797274489383016717182991353216456e-01}
 \ . \label{eq.eomg2}
\end{align}

\begin{figure}[ht]
\resizebox{\columnwidth}{!}{\input{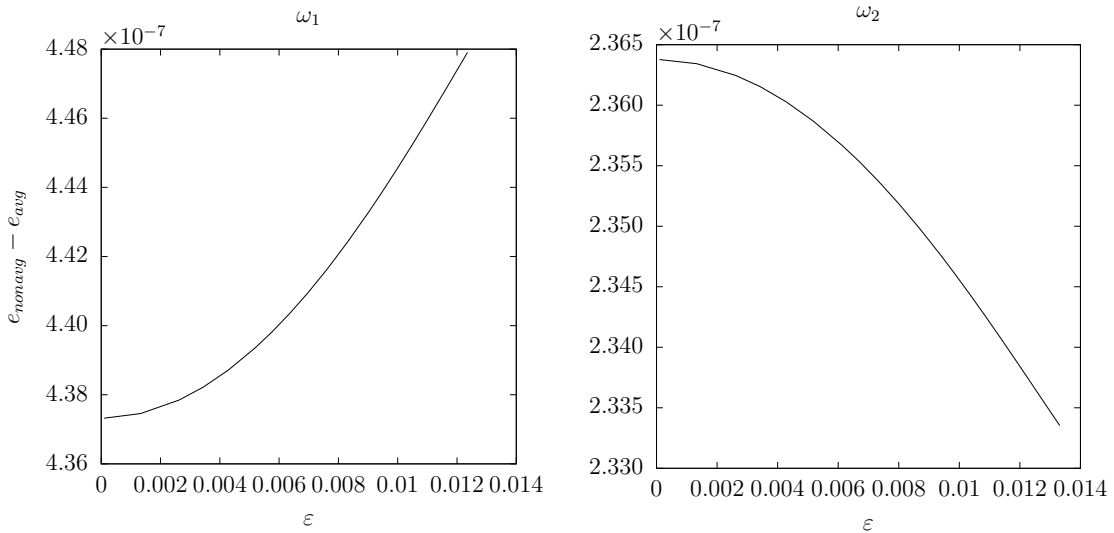}}
 \caption{Difference of the eccentricities for the averaged and
   non-averaged spin-orbit problems \eqref{eq.avg-spinxy} and
   \eqref{eq.spinxy}, respectively, and for the frequencies
   \eqref{eq.omg1} (left panel) and \eqref{eq.omg2} (right panel).}
 \label{fig.avg-navg}
\end{figure}

The critical perturbative and drift parameters,
\eqref{eq.epsomg1}--\eqref{eq.eomg2}, do not substantially differ for
the averaged and non-averaged models (respectively,
\eqref{eq.avg-spinxy} and \eqref{eq.spinxy}). For example,
Figure~\ref{fig.avg-navg} shows that the difference in the drift
parameter (the eccentricity) of the KAM tori, obtained using the
continuation procedure, for the two models is of the order of
$10^{-7}$. Therefore, in single-precision arithmetic, these averaged
and non-averaged drifts may look indistinguishable.

\begin{figure}[ht]
\resizebox{\columnwidth}{!}{\input{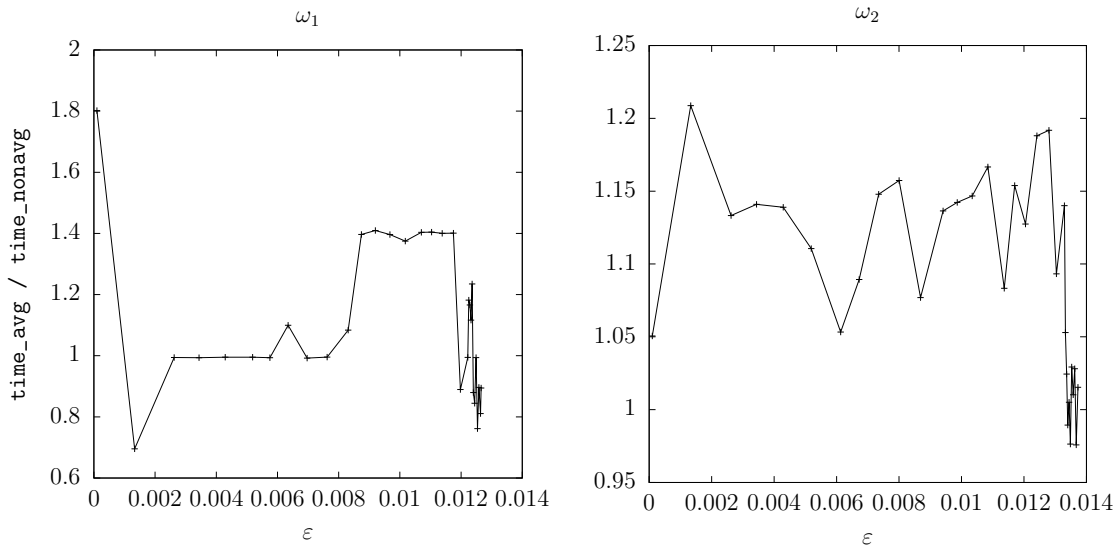}}
 \caption{Time ratios of each successful continuation step for the
   averaged and non-averaged spin-orbit problems \eqref{eq.avg-spinxy}
   and \eqref{eq.spinxy}, respectively, and for the frequencies
   \eqref{eq.omg1} (left panel) and \eqref{eq.omg2} (right panel). It
   run with 32 CPUs.}
 \label{fig.stats-avg-navg}
\end{figure}

The mesh size is the same in the two models and for the studied
frequencies, i.e., Figure~\ref{fig.stats-avg-navg-mesh} is the same
for the averaged model. Also the computational times to compute the
tori are similar. In general, the averaged version is faster far from
the breakdown. Figure~\ref{fig.stats-avg-navg} illustrates the time
ratios (using in all the cases a maximum of 32 CPUs), while
Figure~\ref{fig.stats-avg-navg-mesh} provides the time of each
successful continuation step in $\eps$. We strongly encourage to not
infer general statements from the computational times, since they may differ due to cache sizes, CPU machines, running of other jobs, changes
in the continuation strategies, initial conditions, number of Newton
iterations, accurate control flags, etc.

The main computational bottleneck is on the Newton steps. For
each discretized value of $\theta$ parametrizing the invariant torus,
a numerical integration of the ODE and its variational equations must be performed. Also, different continuation strategies in the stepsize choice and initial guesses can improve the overall time.

Motivated by these remarks, in what follows, we are going to provide
only the plots for the non-averaged spin-orbit problem.  The plots of
the two models are visually identical, since the difference between
them in the regimes we study is about 5 orders of magnitude smaller
than the main effects.

\section{Studies of the  breakdown of tori}
\label{sec:studybreakdown}

Our results on the behavior at breakdown of the spin-orbit map are
described in Sections~\ref{sec:bundles} and \ref{sec:breaksob}, where
we have taken the algorithms to their limits of validity. Our study is
motivated by different reasons:
\begin{enumerate}
\renewcommand{\theenumi}{\sl \roman{enumi}}
\renewcommand{\labelenumi}{(\theenumi)}
  \item the spin-orbit problem provides an excellent testing ground
    for the algorithms, which is far from the academic standard models;
  \item the a-posteriori method allows to compute with confidence very
    close to the breakdown;
  \item the phenomena at the breakdown of validity of KAM theory have
    a great mathematical interest, also because they are related to
    renormalization group theory (\cite{OstRanSetSig-83,
      RanOstSet-82}).

For us, the phenomena at breakdown is a side issue, since the main
goal of this paper is to develop a methodology to maintain accuracy
and reliability even up to breakdown in a model of astronomical
interest. Clearly, the phenomena happening at the breakdown are of certain interest and deserve another study.
\end{enumerate}
We anticipate that the results that we obtain are that the breakdown
of the dissipative spin-orbit problem does not conform to the
behaviors that have been previously found in the literature, \cite{MacKay82,RanOstSet-82}.

\subsection{Some rigorous results on the breakdown}
% There are rigorous results that limit what happens at breakdown and hence their effect on the numerical experiments.

The rigorous result in \cite{CallejaCdlL2013} has a remarkable
consequence concerning the behavior of the torus, and of the stable
and tangent bundles at breakdown.  Precisely, if:
\begin{enumerate}
\renewcommand{\theenumi}{\sl \roman{enumi}}
\renewcommand{\labelenumi}{(\theenumi)}
  \item \label{req1} the invariant attractor is smoothly conjugate to
    a rotation,
  \item \label{req2} the invariant torus has a smooth invariant
    direction which is:
  \begin{itemize}
    \item contractive,
    \item with an angle bounded below from the tangent direction,
  \end{itemize}
\end{enumerate}
then the torus can be continued as a smooth curve and the iterative
method specified in \cite{CallejaCdlL2013} will converge for a small
enough perturbation.

As a consequence of the above results, at breakdown at least one of
the above requirements \eqref{req1} or \eqref{req2} has to fail. In
other words, either the hyperbolicity is lost or there is a loss of
regularity.

Hence, our numerical explorations focus on the behavior of the bundles
(Section~\ref{sec:bundles}), the regularity of the tori
(Section~\ref{sec:domains}), and the regularity in
Sections~\ref{sec:breaksob}. Since previous studies found scaling
relations, we have also studied the possibility of scaling invariant
ratios (Section~\ref{sec:breakdownscaling}).

\medskip

There are some rigorous results that further limit the phenomena that
can happen.

As for the breakdown by loss of hyperbolicity, in the spin-orbit
problem the determinant is the friction, since the dynamics on the
circle is a rotation, the stable exponent is given by the friction, so
the only way that hyperbolicity can be lost is if the stable and
tangent bundles become close at some points.  This phenomenon was
found numerically in many examples (\cite{HLlnum, HLlverge}) and
called {\em bundle collapse}.  In \cite{HLlnum} it was also found that
this phenomenon happens often and that it satisfies remarkable scaling
relations. The paper \cite{HLlnum} contains a proof that this
phenomenon indeed happens in some examples. In
\cite{BjerklovS08,FiguerasT20} there is a proof of the scaling
relations of bundle collapse in some models.

As for the breakdown by loss of regularity, we note:
\begin{itemize}
\item If the torus and the bundles are $C^1$, using that the normal
  exponent has to be the fraction, the results in \cite{Fenichel74}
  show that the torus and the bundle are $C^r$ for any $r \in
  \mathbb{N}$.
\item Using \cite{Herman79, Yoccoz84, KhaninS86}, we obtain that the
  dynamics is $C^{r-a}$ conjugate to a rotation for some $a \in \R _+$
  related to the Diophantine exponent.
\item If the torus is $C ^{r-a}$ conjugated to a rotation, then the
  bundles are $C ^{r-2a}$.
\item In \cite{CallejaCdlL2013} it is shown that if the torus and the
  bundles are $C^r$ for $r$ sufficiently high, then the torus and the
  bundles are analytic.
\end{itemize}
Therefore, at the breakdown by loss of regularity, then the regularity
of the torus or the bundle has to go from analytic to less than $C^1$.

We also note that, as we will see, the bundles can be obtained from
the conjugating function by solving cohomology equations, so that the
only way that bundles can lose regularity is if the regularity of the
torus breaks down.

This drastic drop of regularity justifies that we talk about the
breakdown. In other situations, the regularity seems to decrease
gradually (\cite{YaoL21b}).

We also note that, even after the breakdown, there could be other
invariant objects of lower regularity. For example, \cite{CapinskiK20}
presents a topological extension of normal hyperbolicity beyond $C^1$
regularity (the objects produced are continua, not manifolds).  For
dissipative systems, there are topological
(\cite{LeCalvez86,LeCalvez88}) or variational arguments
(\cite{SorrentinoM17}), showing existence of quasi-periodic orbits of
a rotation number, even if they are not dense on a circle.

The closest precedents to our study are \cite{HLlnum, HLlverge} and,
specially \cite{CallejaF2012}, which studies conformally symplectic
systems. In these papers, it was found at the same time a bundle
collapse and a loss of regularity as well as remarkable scaling
relations.

\section{Behavior of invariant bundles}
\label{sec:bundles}
In this section, we study the behavior at breakdown of the stable and
tangent bundles. To this end, let $(K,\ecc)$ be an invariant torus
with frequency $\omg$ for the map $P _\ecc$ given in \eqref{SOmap}.
Let $N(\theta)$, $M (\theta)$, $S(\theta)$ be the quantities defined
by
\begin{equation*}
  \begin{split}
N(\theta) &\bydef (DK(\theta)^\top DK(\theta))^{-1}\ , \\ M(\theta)
&\bydef [DK(\theta)\ |\ J^{-1}\circ K(\theta)\ DK(\theta)
  N(\theta)]\ , \\ S(\theta) & \bydef ((DK N)\circ
T_\omega)^\top(\theta) DP_{\ecc} \circ K(\theta) J^{-1}\circ
K(\theta)DK(\theta)N(\theta)\ ,
  \end{split}
\end{equation*}
where $T_\omega$ denotes the shift by $\omega$:
$T_\omega(\theta)=\theta+\omega$.  Then, one can show
(\cite{CallejaCdlL2013}) that $N$, $M$, $S$, and $P_e$ satisfy the
relation
\begin{equation*}
 D P _\ecc \circ K (\theta) M (\theta) =
 M (\theta + \omg)
 \begin{pmatrix}
  1 & S (\theta) \\ 0 & \lambda
 \end{pmatrix}
\end{equation*}
with $\lambda$ defined in \equ{lambda}.  Similarly, the stable and
tangent bundles, say $E ^s(\theta)$ and $E^c(\theta)$ respectively,
must satisfy the reducibility equation (see\cite{CallejaF2012})
\begin{equation*}
 \label{eq.reducibility}
 D P _\ecc \circ K(\theta) W(\theta) =
 W(\theta + \omg)
 \begin{pmatrix}
  1 & 0 \\ 0 & \lambda
 \end{pmatrix}
\end{equation*}
with $W (\theta) =
\begin{pmatrix}
 E^c(\theta) & E ^s(\theta)
\end{pmatrix}$.

To reduce the cocycle, we introduce the change of variables
\begin{equation*}
 \widehat W (\theta) =
 \begin{pmatrix}
  1 & B (\theta) \\ 0 & 1
 \end{pmatrix}
\end{equation*}
satisfying
\begin{equation*}
 \begin{pmatrix}
  1 & S (\theta) \\ 0 & \lambda
 \end{pmatrix} \widehat W (\theta) = \widehat W (\theta + \omega)
 \begin{pmatrix}
  1 & 0 \\ 0 & \lambda
 \end{pmatrix}.
\end{equation*}
Hence, the unknown function $B(\theta)$ verifies the following
cohomology equation:
\begin{equation}
 \label{eq.cohom3}
 B(\theta) - \lambda B (\theta + \omg) = - S (\theta)\ ,
\end{equation}
which can be solved by expanding in Fourier series whenever $| \lambda
| \ne 1$.  Finally, the tangent and stable bundles can be recovered
from the relation $W(\theta) = M (\theta) \widehat W (\theta)$.
Explicitly, if $K $ has components $(K _1, K _2)$, then
\begin{equation*}
 W(\theta) =
 \begin{pmatrix}
  D K _1(\theta) & D K _1(\theta) B (\theta) - D K _2(\theta) N (\theta) \\
  D K _2(\theta) & D K _2(\theta) B (\theta) + D K _1(\theta) N (\theta)
 \end{pmatrix}.
\end{equation*}
A straightforward computation of the inner products between the
bundles leads to
\begin{equation*}
 \begin{split}
  E^c(\theta)^\top E^s(\theta) &= N(\theta) ^{-1} B(\theta), \\
  E^c(\theta)^\top E^c(\theta) &= N(\theta) ^{-1} , \\
  E^s(\theta)^\top E^s(\theta) &= N(\theta) ^{-1} B(\theta)^2 + N(\theta)\ .
 \end{split}
\end{equation*}
Therefore, if $\alpha(\theta)$ is the angle between the invariant
stable and tangent bundles at the point $\theta \in \T$, then
\begin{equation*}
 \cos \alpha (\theta) = \frac{B(\theta)}{\sqrt{N(\theta)^2 +
     B(\theta)^2}}
\end{equation*}
     or equivalently
\begin{equation}
 \label{eq.tanangle}
 \tan \alpha (\theta) = \frac{N(\theta)}{B(\theta)}\ .
\end{equation}

\begin{figure}[ht]
\resizebox{.99\columnwidth}{!}{\input{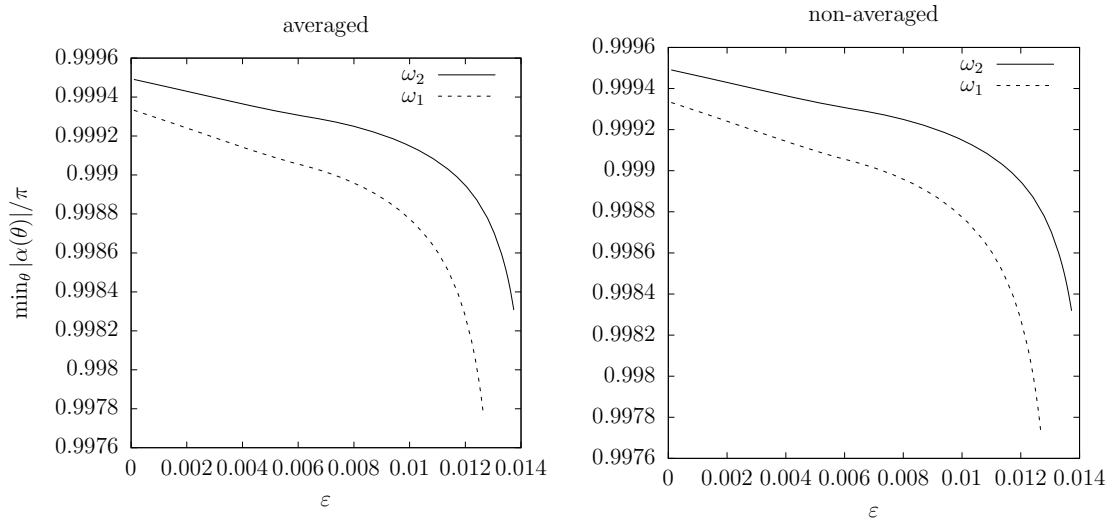}}
 \caption{Minimum angle between the stable and tangent bundles of the
   KAM tori of \equ{A1} (left panel) and \eqref{eq.spinxy} (right
   panel) with $\dis=10^{-3}$, $\ecc$ as drift, and $\eps$ as
   continuation parameter.}
 \label{fig.angle-bundles}
\end{figure}

Figure~\ref{fig.angle-bundles} shows the behavior of the angle
$\alpha$ for the averaged and non-averaged models \eqref{A1} and
\eqref{eq.spinxy}. We notice that the minimum value of the angle of
separation does not go to zero as the parameter $\eps$ approaches the
breakdown threshold.

In \cite{CallejaF2012} it is shown that the angle between the tangent
and stable bundles of the dissipative standard map collided in a very
particular manner, while we have realized that in the spin-orbit
problem the approach to breakdown occurs in a different way.  More
precisely, let us consider the angles w.r.t. the semi-axis $\{x>0\}$,
which are obtained by computing the arc tangent between the second
coordinate and the first one in each of the
bundles. Figure~\ref{fig.bundles} shows these angles for the center
and stable bundles as well as their difference for the $\eps$ values
given in \eqref{eq.epsomg1} and \eqref{eq.epsomg2}. We observe that
the difference between these angles does not coincide as in the
dissipative standard map, yielding a non-collapse bundle near the
breakdown in the sense of \cite{CallejaF2012}. Moreover, the minimum value happens in the same value around $\theta\approx 0.5$ which is coherent with the results in \cite{BjerklovS08}.

\begin{figure}[ht]
\resizebox{.99\columnwidth}{!}{\input{bundles-and-diff.tex}}
 \caption{Angles between the center and stable bundles w.r.t. the
   semi-axis $\{x>0\}$ and its difference for $\ecc$ (\equ{eq.eomg1}
   upper left, \equ{eq.eomg2} lower left) and $\eps$ \equ{eq.epsomg1}
   upper right, \equ{eq.epsomg2} lower right) values corresponding to
   each of the frequencies \eqref{eq.omg1} and \eqref{eq.omg2} of the
   KAM tori of the non-averaged spin-orbit \eqref{eq.spinxy} with
   fixed $\dis=10^{-3}$.}
 \label{fig.bundles}
\end{figure}

\begin{figure}[ht]
\resizebox{.49\columnwidth}{!}{\input{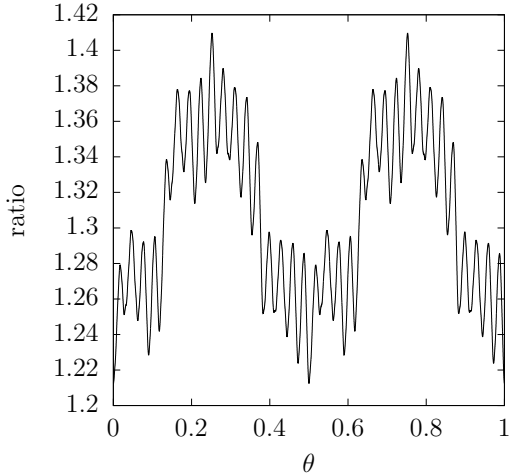}}
 \caption{Quotient between the difference in the angle of the center
   and stable bundles for frequencies \eqref{eq.omg1} and
   \eqref{eq.omg2} respectively, see Figure~\ref{fig.bundles}. As in
   that figure, the computations refer to KAM tori of the non-averaged
   spin-orbit \eqref{eq.spinxy} with fixed $\dis=10^{-3}$.}
 \label{fig.bundles-diff-ratios}
\end{figure}

In conclusion, the separation of the bundles remains bounded away from zero as
we approach the breakdown. On the other hand, the regularity seems to
decrease.  This indicates that the destruction of the torus does not
happen because of destruction of the hyperbolicity, but just because
of loss of regularity.

\subsection{Similarity of the bundle angles for different rotation numbers}

One striking feature of Figure~\ref{fig.bundles} is the similarity of
the position of the bundles for the two frequencies at breakdown.  Of
course, before breakdown they are completely different.  This could be
explained because in renormalization, many phenomena in the fine scale
at breakdown depend only on the tail of the continued fraction of the
expansion. To illustrate this similarity,
Figure~\ref{fig.bundles-diff-ratios} shows the difference in the angle
of center and stable bundles reported in Figure~\ref{fig.bundles}.

The paper \cite{MacKay82} includes a renormalization procedure that
involves not only the mapping, but also the rotation number.  Roughly,
the basic idea of the transformation is to change the map by an
iterate (related to the return map of the rotation sought), and scale
the space and the parameter to breakdown. At the same time, one
changes the required rotation number by applying to it the Gauss map
$\omega \mapsto 1/\omega - [ 1/\omega]$, or equivalently, shifting the
continued fraction expansion. If this renormalization map has an
asymptotic behavior, which is the same for several maps/numbers (even
if the dynamics of the renormalization is more complicated than a
fixed point), then the properties close to the breakdown would only
depend on the tail of the continued fraction expansion. This is
consistent with our data.

\section{Width of analyticity domains.}
\label{sec:domains}

Since the functions we consider are analytic up to the breakdown, a
way to ascertain the breakdown is to explore when the analyticity
domain goes to zero.

A widely applicable method to approximate the analiticity domain of a
periodic function $f$ is to perform a linear regression of the $\log
_{10} |\hat f_k| $ values (one has to discard the first terms).

In our problem, applying the above method, we obtain $\delta_\eps$ for
each value $\eps$ in the continuation. The function $f_\eps$ is
analytic in $-\delta _\eps \log (10)/(2\pi) $.

\begin{figure}[ht]
\resizebox{.99\columnwidth}{!}{\input{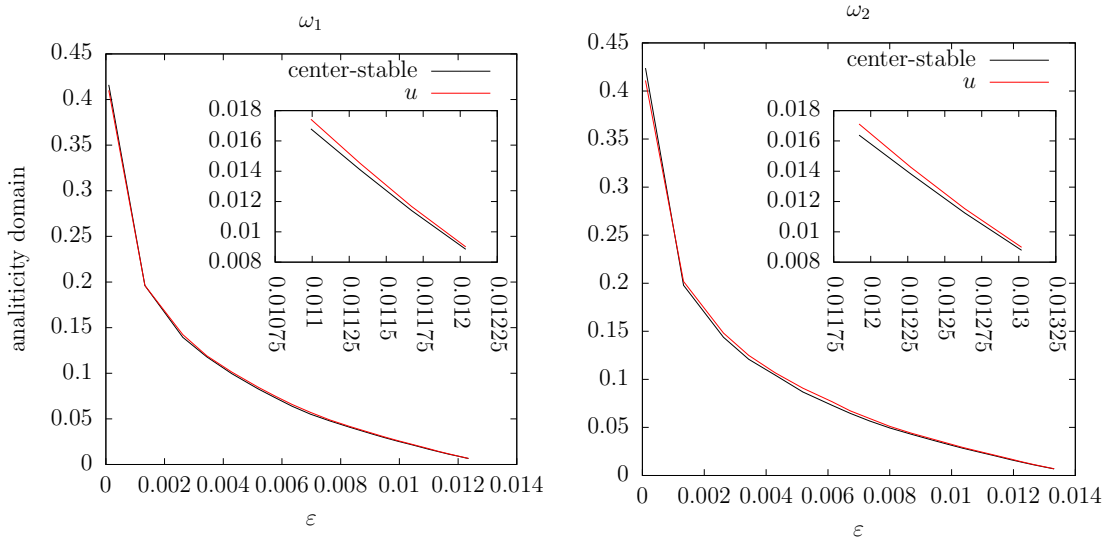}}
 \caption{Maximum size of the analiticity domains of the center minus
   stable angles in Figure~\ref{fig.bundles} and $u$ given in
   \eqref{uu} for each of the frequencies \eqref{eq.omg1} (left panel)
   and \eqref{eq.omg2} (right panel) of the model \eqref{eq.spinxy}
   with fixed $\dis=10^{-3}$, and $\ecc$ as drift parameter. The small
   plots are just zooms in.}
 \label{fig.analiticity}
\end{figure}

Figure~\ref{fig.analiticity} shows the analiticity domain for the
difference of the angles in Figure~\ref{fig.bundles} and also for the
function $u=u(\theta)$ defined as \beq{uu} u(\theta) \bydef K
_1(\theta) -\theta\ , \eeq where $K$ is the spin-orbit embedding of
the torus and $K _1$ its first component.

In both cases, the values tend to zero as $\eps$ approaches the
threshold value.

We note that, since the domain of analyticity has dimensions of
length, if there are scaling relations, then the scalings should go to
zero like a power of the distance of $\eps$ to the breakdown value
\cite{Llave92b}. Note that Figure~\ref{fig.analiticity} does not
support the existence of a scaling behavior.

\section{Sobolev's criterion}
\label{sec:breaksob}
The Sobolev criterion was introduced in \cite{CallejaL09} for
symplectic mappings and extended in \cite{CallejaCdlL2013} to
conformally symplectic mappings.

The rigorous underpinning is very strong. The a-posteriori theorem
(\cite{CallejaL09}) gives an algorithm for the computation of the
analiticity breakdown, based on the fact that, if the approximate
solution is well behaved, there is a true solution near the
approximate one.

We start by giving the definition of Sobolev seminorm as follows.

\begin{defn}\label{defn.sobnorm}
 Let $f$ be a map in $L ^2(\T)$, whose Fourier expansion is written as
 $f(\theta) = \sum _{k \in \Z} \hat f _k e^{2\pi \I k \theta}$, and
 let $\| f \| _{L ^2} \bydef \bigl( \sum _{k \in \Z} | \hat f _k |
 ^2\bigr) ^{1/2}$. Then, if $r$ is a real number, the $r$-th Sobolev
 seminorm is defined by
 \begin{equation}\label{seminorm}
   \| f \| _r \bydef \| \partial _{\theta}^r f \| _{L ^2} =
   \biggl( \sum _{k \in \Z} |2 \pi k| ^{2r} | \hat f _k | ^{2} \biggr) ^{1/2} \
   ,
 \end{equation}
 where $\partial _\theta ^r$ denotes the $r$-th derivative w.r.t. to
 $\theta$ and we have used the Parseval identity to express the $L^2$
 norm of the derivative in terms of Fourier coefficients.
\end{defn}

Note that
\begin{enumerate*}[label={\sl (\roman*)}]
 \item $r$ does not need to be an integer,
 \item the Definition~\ref{defn.sobnorm} is given in terms of a
   seminorm on a space of (smooth) periodic functions, and
 \item a map $f$ in Definition~\ref{defn.sobnorm} will numerically be
   represented by a finite sum, say $f^{\leq L}$ where $L$ denotes the
   meshsize, see Figure~\ref{fig.stats-avg-navg-mesh}.
\end{enumerate*}
Therefore, we compute the seminorm of the truncation of the function
$f$ up to an integer $L$ as
\beq{SN}
 \| f ^{\leq L} \| _r = \biggl( \sum _{|k| \leq L} (2 \pi k) ^{2r} |
 \hat f _k | ^{2} \biggr) ^{1/2} \ .
\eeq

\begin{rem}
  Note that the computation of \eqref{SN} could be very sensitive to
  round-off error, since it involves summing terms of different
  sizes. It is well known that summing first the larger terms and then
  the smaller ones is very affected by round off errors.

  If the floating point of the computer satisfies some elementary
  properties (indeed holding in several popular packages such as
  {\tt MPFR}), then there is a remarkable algorithm due to W. Kahan,
  which computes the sums (in a few extra operations) without round
  off error (\cite[p. 244]{Knuth97}).

  We have made sure that our calculations of the Sobolev seminorm are
  done using this algorithm (it is implemented by default in some
  linear algebra packages).
\end{rem}

%%%%% The rigorous justification of the Sobolev criterion is very strong:

At the basis of the Sobolev criterion there is the following
remark. The KAM theorem developed in \cite{CallejaCdlL2013} for
conformally symplectic systems shows that, if we obtain a solution
with a bounded Sobolev seminorm high enough and the bundles are
separated, then one can continue further. The paper
\cite{CallejaCdlL2013} provides explicit rigorous estimates of what is
the meaning of high enough order (in numerical estimates, one can see
that the blow up happens even for orders smaller than what the
rigorous results ensure).

Therefore, when the parameters approach the threshold, either all the
Sobolev seminorms of the conjugacy of sufficiently high order blow up
or -- inclusive or -- we have that the distance between the stable and
unstable bundles goes to zero.

It is important to remark that the method and its justification remain
valid for any number of dimensions as well as for other problems such
as \emph{nontwist circles} (\cite{GonzalezHL22, Cal-Can-Har-21}).

This method was originally implemented in \cite{CallejaL09} for
symplectic maps, where the angle of the bundles did not enter. It was
found to be very practical, since the continuation algorithms for
quasi-periodic tori do not require adjustment and can run
unattended. One can run different paths of continuation in different
cores and obtain domains of two parameters easily.  Independent
implementations confirming and extending the results to asymmetric
mappings appeared in \cite{FoxM16} (note that the asymmetric maps
considered in \cite{FoxM14} do not have symmetry lines, so that the
methods based on periodic orbits such as Greene's method have great
difficulty).  Other independent implementations for twist maps appear
in \cite{Flesher21}, which explored even a 3 dimensional parameter
space. Implementations for the breakdown of two dimensional tori
appear in \cite{BlassL13, FoxM16}.

\medskip

\begin{figure}[ht]
\resizebox{\columnwidth}{!}{\input{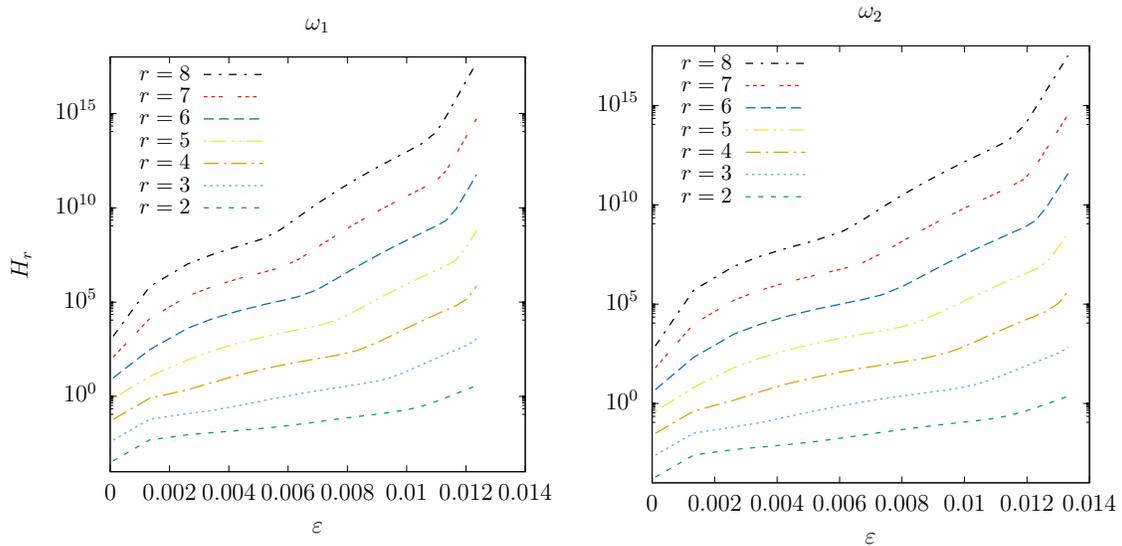}}
 \caption{Sobolev seminorms as in \equ{SN} of the spin-orbit embedding
   $u$ in \equ{uu} for different values $r$ for the frequencies
   \eqref{eq.omg1} (left panel) and \eqref{eq.omg2} (right panel) of
   the non-averaged model \eqref{eq.spinxy}. }
 \label{fig.sob}
\end{figure}

In numerical computations, we consider the function $u$ in \equ{uu}
and we approximate it as $u^{\leq L}$ for a suitable integer $L$.  The
truncated seminorm of $u^{\leq L}$ will be denoted as $H _r \bydef \|
u ^{\leq L} \| _r$. The function $u$ depends on the frequency of the
torus.  Figure~\ref{fig.sob} provides the seminorm values for
different indexes $r$ and for each of the two frequencies considered
in this paper, namely $\omg _1$ in \eqref{eq.omg1} and $\omg _2$ in
\eqref{eq.omg2}. As $\eps$ increases, the seminorms blow-up,
indicating a breakdown when approaching the critical $\eps _c$ values.

\subsection{On other methods to compute the breakdown values}
\label{sec:other}
There have been several methods proposed to study numerically the
breakdown of invariant circles, also for higher dimensions \cite{BustamanteCL23}.  A survey for symplectic systems
appears in \cite{CallejaL10}; some of the methods exclude circles of a
given rotation number, some of them exclude all circles.  An
interesting question is to study which of these symplectic methods can
be adapted to conformally symplectic systems. This adaptation is not
trivial since the dissipative systems require to adjust parameters to
maintain the frequency.  Besides the conceptual problems of
formulation and justification, there are of course practical problems.

We have considered adapting the Greene's criterion \cite{Greene1979}
and the obstruction criterion \cite{OlveraS86}.  The Greene's
criterion requires some adaptation (e.g. redefining the residue,
clarifying the role of parameters (these extensions were  accomplished in \cite{CCFL14}), while the
obstruction criterion is topological and does not need adaptation.
Note also that the obstruction criterion is carefully justified,
while the Greene's criterion has only partial justifications.

Both criteria depend on a detailed study of periodic orbits, whose
rotation number approximates the frequency of the torus. Finding
periodic orbits is not easy, even in standard-like symplectic maps, in
particular when they contain several harmonics (\cite{LomeliC,
  FalcoliniL92b}). For dissipative mappings, Greene's criterion has
been successfully implemented in maps with few harmonics (see
\cite{CallejaC10}, which also includes a comparison with Sobolev
criterion in those cases).

Certainly, the spin-orbit map contains many harmonics, there are no
symmetry lines and it becomes very hard to track periodic orbits. Note
also that one has to adjust the drift parameter; periodic orbits --
with different characteristics -- exist for drifts in intervals
(Arnold tongues).

The obstruction criterion requires a careful computation of the stable
and unstable manifolds of periodic orbits. We have given some effort
to the implementation, but found very hard to follow periodic orbits
of high period in the spin-orbit problem.

At the moment, we are not aware of any computational alternative to
the Sobolev criterion, which is, moreover, rigorously justified.

The Sobolev criterion has the weak point (shared with Greene's
criterion) that the breakdown involves a limit. Hence, a finite
calculation cannot give a rigorous conclusion (calculations involving
extrapolations are standard in numerical analysis, though).  The obstruction
criterion produces rigorous results of non-existence with a finite
calculation (see \cite{CellettiM07} on the non existence of tori of
all rotation numbers). Our KAM theorems produce definite results of
existence with a finite calculation (\cite{CCGL20b}).

\section{Fine properties of the breakdown: scaling relations, renormalization
  group}
\label{sec:breakdownscaling}
One of the consequences of our calculations (see
Section~\ref{sec:breakdown}) is that we can explore the phenomena that
happen at breakdown of KAM theory, which is an area full of open
mathematical problems.

As a consequence of our detailed calculations, we find that the
breakdown of invariant circles in the spin-orbit problem does not fit
with the customary scaling theory, documented in the breakdown of
other mappings \cite{RanOstSet-82,OstRanSetSig-83}.

The emphasis of this paper is in the accurate computation of tori in a
concrete spin-orbit problem, the behavior at breakdown is a (welcome
but unexpected) byproduct.  Therefore, in this paper, we will not
investigate systematically other properties such as universality
(i.e., the behavior of several maps) which are of great importance to
the renormalization community, but which require a point of view
different from the one considered here.

\subsection{Scaling theory and renormalization group}

In this section, we present a very quick overview of scaling theory
and renormalization group. The presentation is informal and we omit
many important considerations such as domains of definitions, scaling
properties of angles, etc.

A very important discovery (\cite{RanOstSet-82,Rand92a, Rand92b,
  Rand92c, MacKay82, Kadanoff, Shenker, FeigenbaumKS82}) in the 1980's
was that, for the systems studied, the breakdown of KAM tori had
similarities with the theory of phase transitions in statistical
mechanics: it satisfies scaling relations with \emph{universal}
exponents (i.e., the exponents appearing in the scaling relations are
the same for similar systems).  These scaling properties and the
universality of the exponents are explained and predicted by
properties of the dynamics of a renormalization group (a dynamical
system in an infinite dimensional space of maps).

\def\teps{{\tilde \eps}}
\def\th{\theta}

If we have a breakdown of a KAM torus for $\eps = \eps_c$, let us
denote by $\teps \bydef \eps - \eps_c $.  Scaling theory tells us that
scaling the parameter is asymptotically the same as scaling the
conjugacy. More precisely, for some numbers $\delta$ and $\eta$, we
have
\begin{equation}\label{scaling1}
  K_{\delta^{-1} \teps}(\th)   \approx K_\teps (\eta \th)
\end{equation}
and there is a limit function $K^*$, such that
\begin{equation}\label{scaling2}
    K_{\delta^{-n} \teps}(\eta^n \th) \rightarrow K^*(\th), \qquad
    \text{as } n \to \infty\ .
\end{equation}
Furthermore, the scaling factors $\delta$ and $\eta$ are independent
of the family considered; besides, the $K^*$ that can appear
considering different families (or taking different initial $\eps$)
belong to a 1-dimensional family of mappings.

This scaling behavior can be explained by assuming that a
renormalization acting on mappings has a non-trivial fixed point with
a 1-dimensional unstable manifold and a codimension 1 stable
manifold. The main object of investigation in our paper is the scaling
behavior, which is a consequence of dynamical assumptions of a
renormalization operator.
As a matter of fact, \eqref{scaling1} becomes exact in the unstable manifold of
renormalization.

We omit the details of the well known derivation of scaling behaviors
from the renormalization picture, since they are not relevant for our
discussion.  We just recall that renormalization theory introduces an
operator in the spaces of maps, which is just a change of scale in
time and in space.  Assuming that there is a fixed point of this
operator, the scaling factor $\delta$ is the unstable eigenvalue at
the fixed point, the scaling $\eta$ is given by properties of the
fixed point and the set of scaling limits is the unstable
manifold. The different maps at the threshold for many families form
the stable manifold of the fixed point under the renormalization
dynamics.

The renormalization analysis also predicts other scaling relations,
such as the scaling properties of the analyticity domain of the
solutions (see Section~\ref{sec:scale}).

The results reported in Section~\ref{sec:scale} show that the
spin-orbit model does not seem to satisfy the scaling relations
\eqref{scaling1}.

\medskip
As a consequence of our computations, we conclude that the dynamics of
the renormalization group in the neighborhood of the breakdown for the
spin-orbit map is more complicated than the simple dynamics observed
before.

Indeed, Figures~\ref{fig.HH1} and \ref{fig.HH2} are compatible with the
renormalization group having some type of non-trivial recurrence
(e.g. a limit cycle or a homoclinic orbit).
\medskip

The results in this paper are completely compatible with the previous
papers reporting other behaviors. In renormalization theory, the
universality of the exponents holds only on some \emph{universality
domain} (an open domain in the space of maps for which the dynamics of
the renormalization is described by a hyperbolic fixed point). Indeed,
one of the questions raised in \cite{MacKay82,MacKaythesis} is
precisely to explore the boundary of the universality domain.

To obtain results in maps with many harmonics it was crucial to have
the method based on Sobolev criterion, which is very robust.  Other
methods, such as Greene's method, seem to depend very well in having
periodic orbits working in a predictable way even up to close to
breakdown. This seems to be linked to the renormalization having a
simple dynamics \cite{LlaveO}.

Since the main goal of this paper is to develop methods to compute
extremely accurately until breakdown in concrete models (where the
computation is challenging), we have not done a systematic study of
easy to implement maps in a neighborhood.  This study is natural from
the point of view of renormalization group, but it goes in a direction
different from the main goal of this paper.

\subsection{Scale invariant observables}\label{sec:scale}
In this section, we look at the scale invariance properties of the
Sobolev seminorms.  Given the relevance of scalings, it will be
important for us to perform measurements on the approach to breakdown
which are scale invariant.

Following \cite{CallejaL10b}, we observe that, by Parseval identity, if
$ \beta > 0, \eta \in \N $, and $u _\eta(\theta) \bydef u(\eta
\theta)$, then
\[
\begin{split}
    \| \beta^{-1} u _\eta \|^2_r & =
    \beta^{-2} \int_0^1 |D^r u(\eta \th) |^2 \, d\th
      =  \beta^{-2}\int_0^1 \eta^{2r-1}  |(D^r u) (\eta \th) |^2 \, d (\eta \th)\\
     & = \beta^{-2} \eta^{2 r -1}  \int_0^\eta | D^r u(\sigma)|^2 \, d \sigma
     =  \beta^{-2} \eta^{2r} \int_0^1  | D^r u(\sigma)|^2 \, d \sigma  \\
    & =   \beta^{-2} \eta^{2r}   \| u  \|^2_r\ .
\end{split}
\]
Even if the derivation above is mathematically exact only when $\eta$
is a large number, similar scaling relations are true in the leading
order for large $\eta$.

Therefore, if we consider observables which are the products of $N$
Sobolev seminorms, we see that, under scalings they transform as:
\begin{equation}\label{eq.scalingprop}
  \prod _{j=1}^N \|\beta^{-1}  u _\eta \| _{r_j}^{\gamma _j} =
 \eta^{ \sum _{j = 1}^N  \gamma_j r _j}\cdot \left( \beta^{-1} \eta^{-1/2} \right)^{ \sum _{i = 1}^N \gamma _i}\cdot
  \prod _{j = 1}^N\| u \| _{r _j} ^{\gamma _j}\ .
\end{equation}

Therefore, if we chose $\gamma _j$, $r _j$ such that
\begin{equation}\label{eq.rgammacond}
  \sum _{i = 1}^N \gamma _i r _i = 0 \qquad \text{and} \qquad \sum _{i = 1}^N  \gamma _i = 0\ ,
\end{equation}
then, the observables \eqref{eq.scalingprop} are invariant under the scaling no
matter what are the scaling parameters $\eta,\beta$.  We refer to
these observables whose value does not change under scaling as
\emph{scale invariant} observables.

If scaling theory applies, we must see that the scaling invariant
observables  are log-periodic in the parameters in the unstable
manifold.  (Compute a scale invariant observable in the
unstable manifold where \eqref{scaling1} is  exact.  We see that they should
satisfy $\phi(\delta \tilde \eps) = \phi(\tilde \eps)$. For
arbitrary families, the $\lambda$-lemma says that the log-periodicity should
be asymptotic near the breakdown).

One can think of the scaling invariant observables as some coordinates
in the space of mappings. If the standard renormalization group
picture (the breakdown corresponds to the stable manifold of a
non-trivial fixed point) applies, one should see that near the
breakdown all the scale invariant variables accumulate near the values
corresponding to the scaling limits (the one-dimensional unstable
manifold of the renormalization group).

\begin{figure}[ht]
 \resizebox{\columnwidth}{!}{\input{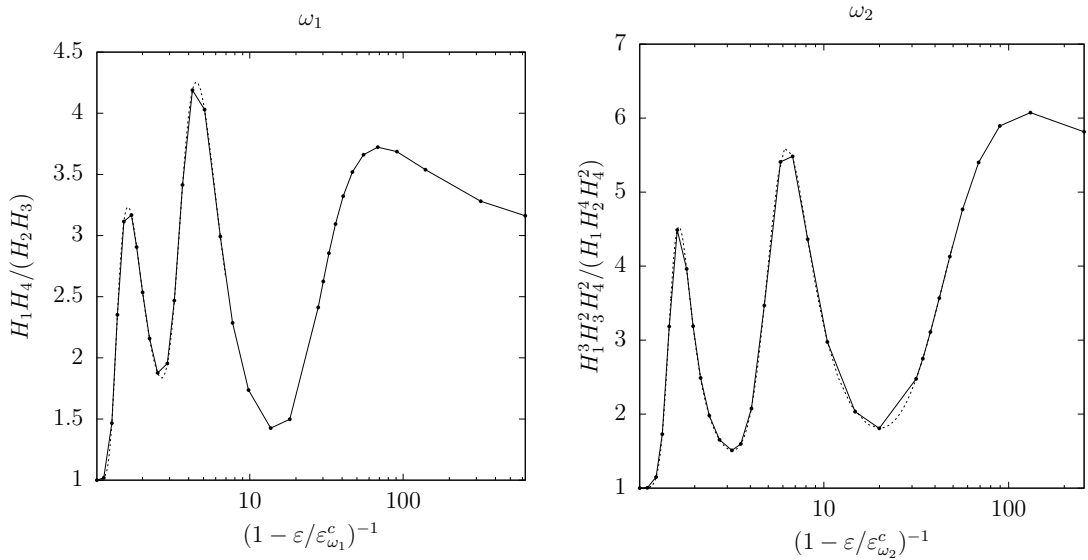}}
 \caption{Scale-invariant ratios of the Sobolev's seminorms $H_r\equiv
   \|u^{\leq L}_r\|$ of the non-averaged spin-orbit problem
   \eqref{eq.spinxy} for $\dis=10^{-3}$ and frequencies
   \eqref{eq.omg1} and \eqref{eq.omg2} versus the continuation
   parameters. The dashed lines represent the same continuation with
   smaller maximum stepsize.}
 \label{fig.HH1}
\end{figure}

\begin{figure}[ht]
\resizebox{\columnwidth}{!}{\input{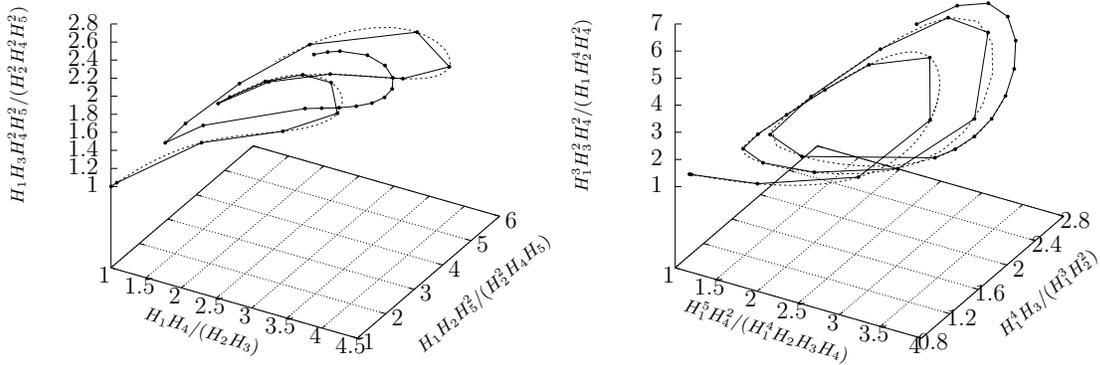}}
 \caption{Scale-invariant ratios of the Sobolev's seminorms
   \eqref{seminorm} for the non-averaged spin-orbit \eqref{eq.spinxy}
   for $\dis=10^{-3}$ and frequencies \eqref{eq.omg1} and
   \eqref{eq.omg2}. Each coordinate line represents a ratio, so that a
   parameter value produces a point in $\R^3$.  Note that the lines
   joining the points are a product of the plotting program.  They
   indicate the progression of the continuation parameter. The dashed
   lines represent the same continuation with smaller maximum
   stepsize.}
 \label{fig.HH2}
\end{figure}

When the breakdown is described by a fixed point of the
renormalization, the ratios of the scaling with $\varrho=0$ tend to a
fixed value monotonically. However, that is not the case for the
spin-orbit model as it is illustrated in Figures~\ref{fig.HH1}
and~\ref{fig.HH2}.

\section{Conclusions} \label{sec:conclusions}
The recent results in \cite{CCGL20a,CCGL20b} show that it is possible
to compute irrational KAM tori in the spin-orbit problem using very
efficient and accurate algorithms. These algorithms are backed up by
rigorous a-posteriori theorems, so that the results are convincing.
Note that the problem is very challenging because it is in the regime
of weak dissipation, which is a singular perturbation of the
conservative spin-orbit problem.

The challenge faced in this work is to take the algorithm to the
limits of validity, and see how one can maintain the reliability and
explore the phenomena that happen.

The main conclusion is that the algorithm maintains the extremely high
accuracy and reliability very close to breakdown.

The calculations of the breakdown in this paper have been based on the
Sobolev criterion of \cite{CallejaL10} for symplectic mappings
(adapted in \cite{CallejaCdlL2013} to conformally symplectic
mappings). We have found it easy to implement and reliable.

One of the results of our exploration is that the breakdown happens in
ways that are qualitatively different to previous explorations.
Notably, we do not observe scaling behaviors and the stable and
tangent directions of the tori remain separated. The main reason why
the tori disappear is just the loss of regularity.

The different behaviors at breakdown can be explained by postulating
that the renormalization map has different dynamics in different
regions in the space of maps. This paper includes tools such as accurate computation and diagnostics like scale invariants that open the possibility of a
more systematic study of the dynamics of renormalization, but this
will require studying breakdown in many families of maps even if they
are not relevant for astrodynamics.

The (very interesting) goal of exploring the renormalization theory
will be postponed, since it is different from the main goal of this
paper: maintaining high precision and reliability in astronomically
motivated models even very close to breakdown.

\subsection*{Conflict of interests}
The authors declare no conflict of interest.

\subsection*{Data availability}
All the data in this manuscript have been generated by the authors. Under reasonable requests, the authors will kindly provide the data from the figures besides the ones provided in the ArXiV version.

% \bibliographystyle{alphaabbr}
%\bibliographystyle{alpha}
%\bibliography{ccglref}

\newcommand{\etalchar}[1]{$^{#1}$}
\def\cprime{$'$} \def\cprime{$'$} \def\cprime{$'$} \def\cprime{$'$}

\end{document}